\documentclass[reqno]{amsproc}

\usepackage{amstext,amsmath,amssymb,amsfonts}
\usepackage[latin1]{inputenc}
\usepackage{epsfig}
\usepackage{hyperref}
\usepackage{color}
\usepackage{tikz}
\usepackage[all,arc,dvips,arrow,curve,cmactex]{xy}
\usepackage{cancel}

\usepackage{latexsym}

\newtheorem{lemma}{Lemma}

\newtheorem{definition}{Definition}
\newtheorem{theorem}{Theorem}

\newtheorem{proposition}{Proposition}

\newcommand{\bea}{\begin{eqnarray}}
\newcommand{\eea}{\end{eqnarray}}
\newcommand{\beq}{\begin{equation}}
\newcommand{\eeq}{\end{equation}}
\newcommand{\enn}{\nonumber \end{equation}}

 \newcommand{\und}{\underline}

 \newcommand{\cG}{\mathcal{G}}
 
 \newcommand{\cV}{\mathcal{V}}
 \newcommand{\cE}{\mathcal{E}}
 \newcommand{\cF}{\mathcal{F}}

 \newcommand{\cR}{\mathcal{R}}
 \newcommand{\cP}{\mathcal{P}}

 \newcommand{\wsig}{{\widetilde{\Sigma}}}

\newcommand{\bsig}{\partial{\Sigma}}

 \newcommand{\sset}{\Subset}

 \newcommand{\inter}{{\rm int}}
 \newcommand{\ext}{{\rm ext}}

\title[Embedding HEGs]{Embedding Half-Edge  Graphs 
in Punctured Surfaces %\footnote{ ICMPA-MPA/2017/...} 
}

\author{Remi C. Avohou}
\address[R.C.A.]{
African Institute for Mathematical Sciences, AIMS-Senegal,
Km 2 Route de Joal - Centre IRD, BP1418, Mbour, Senegal,  
and 
International Chair in Mathematical Physics and Applications,
ICMPA-UNESCO Chair, 072BP50, Cotonou, Rep. of Benin}
\email{avohou.r.cocou@aims-senegal.org}

\author{Joseph Ben Geloun}
\address[J.B.G.]{Laboratoire d'Informatique de Paris Nord, LIPN UMR CNRS 7030, 
Institut Galil\'ee, Universit\'e Paris 13, 
99, avenue Jean-Baptiste Cl\'ement, 
93430 Villetaneuse, France, 
and
International Chair in Mathematical Physics and Applications,
ICMPA-UNESCO Chair, 072BP50, Cotonou, Rep. of Benin}
\email{bengeloun@lipn.univ-paris13.fr}

\author{Mahouton N. Hounkonnou}
\address[M.N.H.]{International Chair in Mathematical Physics and Applications,
ICMPA-UNESCO Chair, 072BP50, Cotonou, Rep. of Benin}
\email{norbert.hounkonnou@cipma.uac.bj}

\begin{document}

\maketitle 

\begin{abstract}
It is  known that  graphs cellularly embedded into surfaces are equivalent to ribbon graphs. In this work, we  generalize 
this statement to broader classes of graphs and surfaces. Half-edge graphs 
extend abstract graphs and are useful in quantum field theory in physics.  On the other hand, ribbon graphs with half-edges 
 generalize ribbon graphs and appear in a different type of field theory emanating from matrix models. 
We then give a sense of embeddings of half-edge graphs 
in punctured surfaces and determine (minimal/maximal) conditions for an equivalence between these embeddings and half-edge ribbon graphs. 
 Given some assumptions on the embedding, the geometric dual of a cellularly embedded half-edge graph is also identified. From that point, the duality can be extended
 to half-edge ribbon graphs. Finally, we address correspondences between polynomial invariants evaluated  on dual half-edge ribbon graphs.
\\


\end{abstract}

\tableofcontents

\section{Introduction}

Graphs embedded in surfaces have been studied
in different contexts with applications
ranging from combinatorics, geometry
to computer science (see \cite{joamo} and the reviews \cite{verdi} and \cite{verdi2} for a detailed account on this active subject). 
The correspondence between embeddings of graphs in 
surfaces and ribbon graphs can be  traced back to the work by Heffter in \cite{Hef}.
Cellular embeddings of graphs, as discussed in \cite{hof}\cite{gt}, are particular graph embeddings such that the removal of the graph from the surface
decomposes the surface in spaces homeomorphic to discs.  
An advantage of working with ribbon graphs is that they
form a stable class under usual edge operations such as edge deletion and contraction.
Meanwhile, to remove an edge in an embedded graph might
result in a loss of the cellular decomposition of
the surface. 
It is therefore useful to have several descriptions
of the same object and use its most convenient
characterization according to the context.

Extending abstract graphs, half-edge graphs (HEGs) 
have appeared in modern physics as Feynman graphs of quantum field theories (see \cite{riv} for a review of the subject). 
Different types of field theories generate different
types of HEGs. Depending on the nature of the field (scalar, vector, matrix or tensor valued), more ``exotic'' field theories
have Feynman graphs with a lot more structure than 
abstract graphs. This is the case of matrix models (see, for instance, \cite{dif}) and  noncommutative field theory (consult the review \cite{krf}) with 
Feynman graphs appearing as half-edge ribbon graphs (HERGs), and 
of tensor models with their Feynman graphs as generalized HEGs 
discovered in \cite{ambj} and called stranded graphs
in \cite{avohou} and \cite{avo2}. 

Investigations on HERGs are still active. Formally, HERGs have been studied
using combinatorial maps in \cite{krf}. Among
other results obtained in that work, the partial duality by Chmutov \cite{chmu} was generalized to HERGs. This duality 
was the stepping stone to find in \cite{krf} a Tutte-like polynomial invariant 
for HERGs satisfying a 4-term recurrence relation. 
In a different perspective, HERGs have been also defined 
by gluing of discs along their boundary in \cite{avohou} and a polynomial invariant generalizing Tutte and Bollob\'as-Riordan (BR) polynomials (see \cite{bollo} and \cite{bollo2}) was found in 
the same work.

While it is clear that half-edges have interesting combinatorial properties, one could ask if they can be useful to  topology as well. 
In this paper, we show that half-edges can be used
to encode punctures on a surface.  
To start, we construct HEG cellular embeddings in punctured surfaces using
the so-called regular embedding \cite{bandi}. 
Consider then the underlying HEG of a HERG obtained by keeping its vertex, 
edge and half-edge sets  and the incidence relation between them. 
We show that a HERG can be uniquely associated with a  
cellular embedding of its underlying HEG in a punctured
surface of minimal genus and minimal number of 
punctures. Under some conditions, the number of punctures in 
the surface can be ``maximal," and we can again identify a unique cellular
embedding of the underlying HEG to which a HERG corresponds. 
 That mapping between HERGs and HEG cellular embedding emanates from
the combinatorial distinction between the boundary components of a HERG (seen as a surface with boundary), introduced by half-edges. Theorem \ref{theo} summarizes that main result. 
We then determine conditions
for which the geometric dual of a HEG
cellularly embedded in a punctured surface exists, is essentially 
unique and defines in return a cellular embedding in the same surface.
Theorem \ref{Propo:dual}, another main result, provides the 
construction of the dual of a HEG cellularly embedding. 
Surprisingly, a case when the dual HEG turns out to be 
well defined occurs when the number of punctures in 
the surface is, in the same above sense, maximal.
Based on results on polynomial invariant on HERGs \cite{avohou},
we finally study duality relations between polynomial invariants
evaluated on dual HERGs. We find two nontrivial instances where a mapping
between these polynomials can be performed. 
Theorem \ref{Theo:duality1} and Theorem \ref{theo:dualerg} give 
new relations between  invariant polynomials calculated on dual HERGs,
thereby providing generalizations of a similar relation revealed in \cite{bollo}. 

The paper is organized as follows. The next section 
reviews surfaces, and graph cellular
embeddings in surfaces and  sets up 
of our notations. In section \ref{sect:cellinpunct},
we first review HEGs and then define HEG cellular embeddings 
in punctured surfaces (Definition \ref{def:celheg}). We also list a few consequences
of our definitions.  The paragraph dealing with HERGs contains
a first main result which is Theorem \ref{theo}. 
In section \ref{sect:geomdual}, we construct the geometric dual of a HEG cellular
embedding in the same punctured surface and Theorems 
\ref{Propo:dual} and \ref{theo:dual} are main results concerning this analysis. 
Finally, in section \ref{sect:dual}, we investigate relationships
between polynomial invariants on dual HERGs. 
We identify two situations where this relationship can be made
explicit.

\section{Surfaces, graph cellular embeddings, ribbon graphs}
\label{sect:surfcell}

In this section, we first review closed and punctured surfaces and their equivalence up to homeomorphism. 
Setting up also our notations, we then  
quickly address cellular embeddings of graphs 
in surfaces and their relationship with 
ribbon graphs. 

\medskip 

\noindent{\bf Surfaces -}
Let $\Sigma$ be a closed connected compact surface of genus $g(\Sigma)$ and $\chi (\Sigma)$ its Euler characteristic. We have
\bea
\chi (\Sigma) = \left\{\begin{array}{ll} 
2-2g(\Sigma) &  {\text{if}}\,\, \Sigma \,\, {\text{is orientable}},\\
2-g(\Sigma)&  {\text{if}}\,\, \Sigma \,\, {\text{is non-orientable}}.\\
\end{array} \right. 
\eea

Let $\Sigma_1$ and $\Sigma_2$ be closed connected compact surfaces. Then $\Sigma_1$ and
$\Sigma_2$ are homeomorphic if and only if they are both orientable or both non-orientable,
and they have the same genus.

A punctured surface $\Sigma$ is a surface obtained after removing a finite number of closed discs (equivalently, up to homotopy,  a finite number of points) in a closed surface. Each boundary component of $\Sigma$ is homeomorphic to a circle and, by capping off the punctures that 
is inserting back the closed discs in $\Sigma$, we obtain a closed surface denoted $\wsig$. We will only be interested in the case
of surfaces $\Sigma$ yielding after capping off a surface $\wsig$ which is compact. 
The  boundary of $\Sigma$ is denoted $\bsig$. 
The genus of $\Sigma$ or of $\Sigma\cup\bsig$  is defined to be the genus of $\wsig$.

Let $\Sigma_1$ and $\Sigma_2$ be connected 
punctured surfaces, and $\wsig_1$ and $\wsig_2$ be the closed connected compact surfaces obtained by capping off the punctures in $\Sigma_1$ and $\Sigma_2$, 
respectively. Then $\Sigma_1$ and $\Sigma_2$  (respectively, $\Sigma_1\cup \bsig_1$ and $\Sigma_2\cup\bsig_2$) are homeomorphic 
if and only if they have the same number of boundary components, and $\wsig_1$ and $\wsig_2$ are homeomorphic.

In the following, a surface will be chosen connected, the general upshot for the non connected
case will be directly inferred from that point.  

\medskip 

\noindent{\bf Graphs -}
A graph $\cG(\cV, \cE)$, or shortly $\cG$, is defined by a vertex set $\cV$, an edge set $\cE$ and an incidence relation between $\cE$ and
$\cV$ (an edge is mapped to a pair of vertices or a vertex in the case 
of a loop). We will first focus on connected graphs and then extend the results to the non connected case.

A graph isomorphism  between $\cG(\cV,\cE)$ and $\cG'(\cV',\cE')$ is a bijection between the vertex sets $\cV$ and $\cV'$, and a bijection 
between the edge sets $\cE$ and $\cE'$ such that two vertices $u$ and $v$ in $\cV$ are adjacent in $\cG$ if and only if their images in $\cV'$ are adjacent in $\cG'$. We denote this isomorphism by $\psi: \cG \rightarrow \cG'$ and say that $\cG$ is equivalent to $\cG'$.

Let $\cG(\cV, \cE)$ be a graph. We associate an underlying topological space $|\cG|$ with the graph $\cG(\cV, \cE)$ as follows (see, for instance \cite{gt,mohar}). Take the sets $\cV$ and $\cE$ each endowed with discrete topology. The space $|\cG|$ is the topological (identification) space $\cV\cup (\cE\times [0,1])$ obtained after the following identification: for each $e\in \cE$, the points $(e, 0)$ and $(e, 1)$ 
are identified with one of the end-vertices of $e$, where for loops, both end-vertices are assumed to coincide. 
 The graph isomorphism $\psi$ extends to graph homeomorphism $\psi_{\rm top}:|\cG|\to |\cG'|$ for topological graphs in a way compatible with the incidence relation. This means that $\psi$ extends to an homeomorphism if each edge $e\times [0,1]$ of $|\cG|$ is homeomorphically mapped to its image $\psi(e)\times [0,1]$ in $\cG'$ with its end vertices or vertex mapped correspondingly.
 The standard textbook of Gross and Tucker, \cite{gt}, gives a survey of topological 
 graph theory. 
Usually we simply write $\cG$ instead of $|\cG|$, providing no confusion arises.

\medskip 

\noindent{\bf Graph cellular embeddings -}
The following definitions of (cellular) embeddings in 
surfaces are withdrawn from \cite{mohar} and \cite{hof}.
2-cells are spaces homeomorphics to open discs. 
For simplicity, we sometimes identify them with open discs. 

\begin{definition}\label{emb}
Let $\cG$ be a graph and $\Sigma$ a surface. An embedding of $\cG$ into $\Sigma$ 
is a continuous map $\phi: \cG \rightarrow \Sigma$ such that the restriction $\phi: \cG \mapsto \phi(\cG)$ is an homeomorphism.
We shall shortly denote both the embedding and the embedded graph by $\cG\subset \Sigma$. 
\end{definition}

\begin{definition}
Let $\cG$ be a graph and $\Sigma$ be a closed connected compact surface. 
An embedding $\cG\subset \Sigma$ is cellular if $\Sigma-\phi(\cG)$  is a disjoint union of $2$-cells.  $\Sigma-\phi(\cG)$ 
is called a cellulation of the surface $\Sigma$
and a 2-cell of the cellulation is called a face of $\cG$. 
\end{definition}

\begin{definition}\label{def:celemb}
We say that two cellularly embedded graphs $\cG_1\subset \Sigma_1$ and $\cG_2\subset \Sigma_2$ are equivalent,  
if there is a homeomorphism $\eta: \Sigma_1 \rightarrow \Sigma_2 $ (which is orientation preserving when $\Sigma_1$ is orientable)
with the property that  $\eta|_{\cG_1}: \cG_1\rightarrow \cG_2 $ is a homeomorphism.
\end{definition}

A cellular embedding $\cG\subset \Sigma$  is said to be orientable if $\Sigma$ is orientable, otherwise we say that $\cG\subset \Sigma$ is non-orientable. If $\cG$ is connected, the
genus, $g(\cG)$, of $\cG$ is the genus of $\Sigma$. A cellularly embedded graph $\cG\subset \Sigma$ is a planar graph if $\Sigma$ is the $2$-sphere.

The Euler characteristic, $\chi(\cG)$, of a cellularly embedded graph $\cG\subset \Sigma$, is defined by
\bea\label{charac}
\chi(\cG) = v(\cG) - e(\cG) + f (\cG),
\eea
where $v(\cG)$, $e(\cG)$, and $f (\cG)$ are respectively the number of vertices, edges and faces of $\cG$.
The Euler characteristic is related to the Euler genus 
by 
\bea\label{gen}
\chi(\cG) = 2 - \gamma(\cG). 
\eea
The above formula extends to a non connected
graph cellularly embedded in a closed connected
compact surface by summing over connected components.
We get:
\bea\label{genk}
\chi(\cG) = 2k(\cG) - \gamma(\cG),
\eea
where $k(\cG)$ is the number of connected components of $\cG$. 

\medskip

\noindent{\bf Ribbon graphs -}
We adopt here the definition by Bollob\'as and Riordan in \cite{bollo} of ribbon graphs. 
 A ribbon graph $G(V,E)$, or simply $G$, is a (not necessarily orientable) surface with boundary represented as the union of two sets of closed topological discs called vertices
 and edges such that vertices and edges intersect by disjoint line segments; each such a line segment lies on the boundary of precisely one vertex and one edge,
and every edge contains exactly two such line segments.
A ribbon graph $G$ naturally has an underlying graph $\cG(V,E)$
that is obtained by keeping only the vertex and edge sets and the incidence between vertices
and edges. We again work with 
connected ribbon graphs and the results will be directly 
extended for non connected ribbon graphs. 

Graphs cellularly embedded in surfaces are equivalent to ribbon graphs and 
this equivalence is established  in the following way. 
To each cellularly embedded graph in a surface, we assign a ribbon graph by taking a  neighborhood strip of the graph in the surface. Reciprocally, given a connected ribbon graph $G$ that we regard as a surface with boundary, 
we cap off that surface by gluing discs along the boundary  components of the ribbon graph. This yields
a closed connected compact surface $\Sigma$ the genus of which, $\gamma(\Sigma)$, is $\gamma(G)$ the genus  of the ribbon graph $G$: 
\bea\label{rgen}
\chi(G) = v(G) - e(G) + bc(G) = 2 - \gamma(G),
\eea
where $\chi(G)$ is the Euler characteristics of $G$ and $v(G)$, $e(G)$, $bc(G)$, are respectively the number of vertices, edges and boundary components of $G$. 
Hence, we have a closed connected compact surface $\Sigma$ endowed already with a cellular decomposition along the underlying graph $\cG$ of $G$. It is also direct to 
observe that
the neighborhood of $\cG$ in $\Sigma$ gives
rise to $G$ again. Noting that the set of faces of $\cG$ is the set of boundary components of $G$, we henceforth call a boundary component a face of $G$, then harmonize our notations and write $bc(G)=f(G)$. 
Importantly, the construction of the cellular
embedding $\cG\subset\Sigma$ stemming from $G$
is \emph{minimal} in the sense that $\Sigma$ is 
the closed connected compact surface with minimum genus in which $\cG$ could be embedded such that the neighborhood $G$ of $\cG$ has the same genus.  

The Euler formula \eqref{gen} generalizes for a ribbon graph $G$ with $k(G)$ connected components as: 
\bea\label{rgenk}
\chi(G) = v(G) - e(G) + f(G) = 2k(G) - \gamma(G). 
\eea

\medskip 

\noindent{\bf Cellular embeddings of  graphs in punctured surfaces -}
We address now the embeddings in a punctured surface of a graph. Note
that the following definition of a cellular embedding of graphs in punctured surfaces differs from 
that of \cite{verdi}. We ensure, for instance,
that the embedding occurs in open surfaces
while, in that work, the boundary is included
in the topological space of the surface.  
In the next definition, 2-cells with punctures
are spaces homeomorphic to open discs 
where we remove some closed discs.  

\begin{definition} 
Let $\cG$ be a graph and $\Sigma$ be a  
 punctured surface. 
An embedding  $\cG\subset \Sigma$ is cellular if $\Sigma-\phi(\cG)$  is a disjoint union of 2-cells possibly with punctures. 
Any 2-cell (with or without punctures) of the cellulation $\Sigma-\phi(\cG)$ is called a face of the embedded graph $\cG$. 
\end{definition}

\begin{figure}[h]
 \centering
     \begin{minipage}[t]{.9\textwidth}
      \centering
\includegraphics[angle=0, width=8cm, height=3cm]{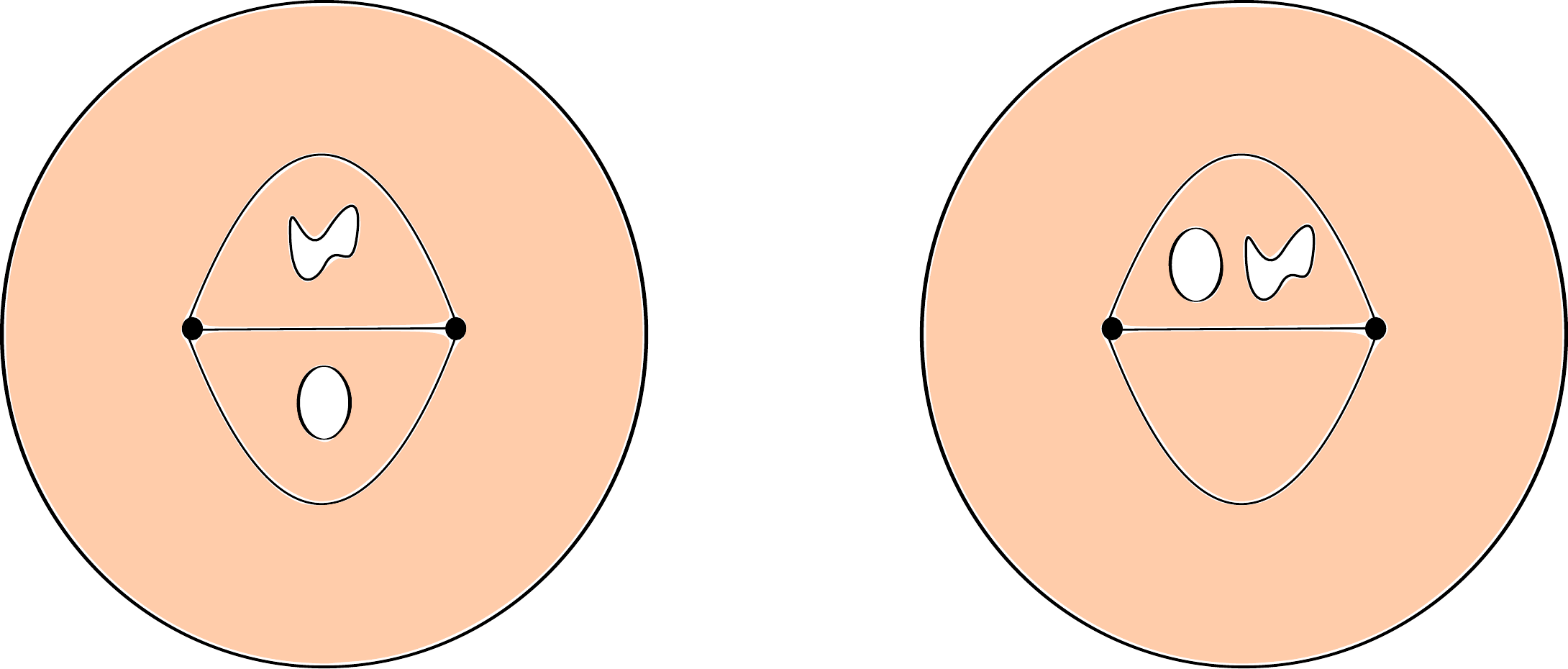}
\caption{ {\small Two equivalent  cellular embeddings
of the same graph in a surface with two punctures (in white). 2-cells may 
contain different number of boundary circles. }}
\label{fig:geps}
\end{minipage}
\end{figure}

Let $\Sigma_1$ and $\Sigma_2$ be two  punctured surfaces. 
Two cellularly embedded graphs $\cG_1\subset \Sigma_1$ and $\cG_2\subset \Sigma_2$ are equivalent if they obey 
Definition \ref{def:celemb}, keeping in mind that surfaces refer now to punctured surfaces. 
One may wonder about the distribution of boundary 
circles in the 2-cells of $\Sigma_1-\phi_1(\cG_1)$
and in the  2-cells of $\Sigma_2-\phi_2(\cG_2)$ which 
may differ (see an example in Figure \ref{fig:geps}).
The above equivalence states that we work 
up to a distribution of boundary
circles in the 2-cells. Said differently, 
a cellular embedding of a graph in a punctured surface $\Sigma$
can be obtained by ``puncturing" after cellular embedding
of a graph in the capping off of $\Sigma$ provided punctures 
are inserted in $\Sigma-\phi(\cG)$.  

Let $\cG\subset \Sigma$ be a graph cellularly embedded in
connected  punctured surface $\Sigma$. $\cG$ is said to be orientable if $\wsig$ is orientable; otherwise we say that $\cG$ is non-orientable. Because $\cG$ is connected, the genus, $g(\cG)$, of $\cG$ is the genus of $\widetilde{\Sigma}$. 
A cellularly embedded graph $\cG\subset \Sigma$ is a plane graph if $\widetilde{\Sigma}$ is the $2$-sphere.
The Euler characteristics and Euler genus for a graph 
embedded in a punctured surface have the 
same formula as  \eqref{charac} and \eqref{gen}
(where $f(\cG)$ now counts the number of all discs including those with punctures) or, in the case of 
a graph with many connected components, as in
\eqref{genk}.

\section{Cellular embeddings of half-edge graphs in punctured
surfaces}
\label{sect:cellinpunct}

We first introduce half-edge graphs  
and then define embeddings of those in punctured surfaces.

\medskip 

\noindent{\bf Half-edge graphs (HEGs) -} We will use notations and conventions of \cite{avohou}. 

\begin{definition}[HEG]
A HEG $\cG(\cV,\cE,h)$, or at times just $\cG_h$, is a graph 
$\cG(\cV,\cE)$, with a set $h$, called the set of half-edges, and a mapping $i:h \to \cV$ called incidence relation which associates each half-edge with a unique vertex. The graph $\cG$ is 
called the underlying graph of $\cG_h$. 
\end{definition}

A HEG isomorphism between $\cG_{h}=\cG(\cV,\cE,h)$ and $\cG'_{h'} =\cG'(\cV',\cE',h')$  is
graph isomorphism $\psi:\cG \to \cG'$ and a bijection between the 
half-edge sets $h$ and $h'$ such that any half-edge $h_0\in h$ is incident to a vertex $v$ in $\cV$ if and only if the corresponding half-edge $h_0'\in h'$ incident to the image of $v$ in $\cV'$.

\begin{figure}[h]
 \centering
     \begin{minipage}[t]{.9\textwidth}
      \centering
\includegraphics[angle=0, width=2cm, height=1.5cm]{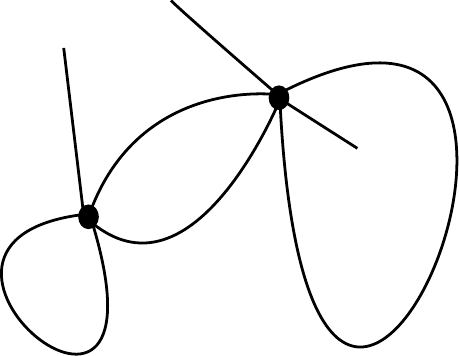}
\caption{ {\small Representation of a HEG with 3 half-edges. }}
\label{fig:heg}
\end{minipage}
\end{figure}

HEGs can be represented in a similar way that abstract
graphs are represented by drawings. To draw 
a HEG, first represent its underlying graph and
then add the set of half-edges represented by a set of segments; each half-edge is incident to a unique vertex without forming a loop.  Figure \ref{fig:heg} illustrates a HEG. 

\begin{definition}[Completed and pruned graphs]\label{completeprun}
Let  $\cG_h=\cG(\cV, \cE,h)$ be an HEG and $h_0\in h$ be one of its 
half-edges incident to a vertex $v_0$. 

 \emph{Completing} $h_0\in h $ in $\cG_h$ is the operation which replaces $h_0$ by an edge $e_0$ by adding a new vertex $v_1$ in the vertex set of $\cG_h$ such that $e_0$ is incident to $v_0$ and $v_1$. 
 The \emph{completed graph} $\bar\cG=\bar\cG(\cV \cup \cV_h, \cE\cup h)$, 
 with $|\cV_h|=|h|$,  
is the graph obtained after  completing all half-edges in $\cG_h$. 

Consider a leaf $v_0\in \cV $ in a HEG $\cG_h$ and $e_0\in \cE$ the edge
incide to $v_0$. 

Pruning $e_0\in \cE$ in $\cG_h$ is the operation which replaces $e_0$ by
a half-edge $h_0$ by removing $v_0$ from the vertex set of $\cG_h$. 
Let $\cV_0$ be a subset of leaves in $\cV$, and $\cE_0=\cE_0(\cV_0)$ be the set of all edges
incident to the leaves in $\cV_0$. 
 The \emph{pruned HEG} with respect to $\cV_0$,  $\und{\cG}_{\,h'(\cV_0)}=\und{\cG}(\cV\setminus \cV_0, \cE\setminus \cE_0, h'(\cV_0) = h\cup \cE_0)$, is the HEG obtained from $\cG_h$ by pruning all leaves in $\cV_0$. 
\end{definition}

Thus completing a half-edge is simply ``promoting'' it as an edge. 
The incidence relation in $\bar\cG$ is an extension of the incidence
relation between edges and vertices in $\cG_h$ by changing $i : h \to \cV$ into $i': h \to \cV \times \cV_h $, such that $i'(h_0)= (i(h_0),v_1)$, and a restriction 
of the incidence relation between half-edges and vertices to an empty mapping. 
Pruning an edge is the inverse operation of completing, that is ``downgrading'' 
an edge as an half-edge. 
 The incidence relation after pruning is a restriction of the
incidence relation between edges and vertices and extension 
of the incidence relation between half-edges and vertices of the former HEG.  
The fact that $\und{\cG}_{\,h'(\cV_0)}$ is a HEG can be then simply verified.

\begin{proposition}
\label{lem:poiujh} 
Let $\cG_h = \cG(\cV,\cE,h)$ be a HEG. 

(1) There is a unique completed graph $\bar\cG$ associated 
with $\cG_h$. 

(2) Fixing $\cV_0$, there is a unique pruned graph $\und{\cG}_{\,h'(V_0)}$ associated
with $\cG_h$. 

(3) The pruned HEG with respect to $\cV_h$ of the completed 
graph $\bar\cG(\cV \cup \cV_h, \cE\cup h)$ is isomorphic to the HEG $\cG_h$. 

\end{proposition}
\proof 
The two first statements are immediate. 
The pruned HEG with respect to $\cV_h$ of 
$\bar\cG(\cV \cup \cV_h, \cE\cup h)$ can be written as
follows 
\bea
\und{(\bar\cG)}_{\,h'(\cV_h)}= \und{(\bar\cG)}
(\cV \cup \cV_h\setminus \cV_h, \cE\cup h \setminus h, h' = \emptyset \cup h)
\,. 
\eea
Thus $\und{(\bar\cG)}_{\,h'(\cV_h)}$ has the same vertex, edge 
and half-edge sets as $\cG_h$. The incidence relation between 
vertices and edges which have been not concerned  by 
the completing and pruning procedures remains
unchanged in $\und{(\bar\cG)}_{\,h'(\cV_h)}$ and $\cG_{h}$. 
The incidence relation between vertices and completed half-edges 
brought by the completing procedure gets restricted by the pruning procedure. 
This implies that the relation between vertices and half-edges remains also the same in 
both  $\und{(\bar\cG)}_{\,h'(\cV_h)}$ and $\cG_h$. 
 
\qed

\begin{figure}[h]
 \centering
     \begin{minipage}[t]{.9\textwidth}
      \centering
\includegraphics[angle=0, width=2cm, height=1.5cm]{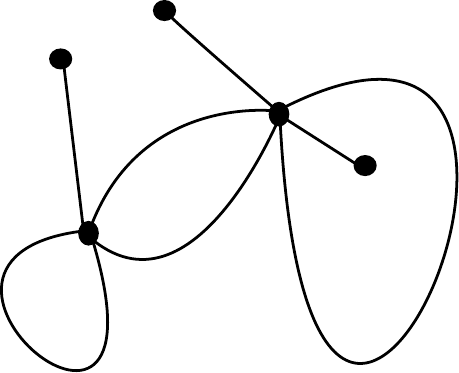}
\caption{ {\small The completed graph of the HEG 
of Figure \ref{fig:heg}: we add new vertices
and turn half-edges to edges. }}
\label{fig:hegstar}
\end{minipage}
\end{figure}

An illustration of the completing of the HEG of Figure \ref{fig:heg} is given  by the graph of  
Figure \ref{fig:hegstar}. Thinking about the inverse operation, i.e to find a HEG $\cG_h$ from a given graph $\cG$,  depending on the number of leaves 
in $\cG$, we can associate a finite number of (possible no) HEGs with $\cG$.

We associate $\bar\cG$ with its underlying topological space $|\bar\cG|$
that we denote again $\bar\cG$, for simplicity. 
Note that we could have introduced
a topology on $\cG_h$ itself, but there is no
need for that since there is now enough data 
to proceed further.

\medskip

\noindent{\bf HEG cellular embeddings -}
 An embedding of a HEG $\cG_h$ in a punctured surface $\Sigma$ follows once again Definition \ref{emb}. It
remains to define the notion of cellulation of 
a punctured surface along a HEG. A way to achieve this
and that further bears interesting consequences is  given by the following. 

 \begin{definition}[$\cV$-regular embedding]
Let $\Sigma$ be a punctured surface with 
 boundary $\partial \Sigma$, such that its capping off 
 gives $\wsig$ a closed connected compact surface. 
Consider a graph $\cG(\cV\cup\cV',\cE)$ with a partition of its vertex set  as shown. 

A $\cV$-regular embedding of $\cG$ in $\Sigma \cup\bsig$ is 
an embedding $\cG \subset \Sigma \cup\bsig$ such that 
  $\wsig-\phi(\cG)$  is a disjoint union of 2-cells
and  $\bsig \cap \phi(\cG)= \phi(\cV)$.  
\end{definition}
 \begin{definition}
 Consider $\bar\cG$ the completed graph of a HEG $\cG_h$
and a punctured surface $\Sigma$. 
  A regular embedding of $\bar\cG $ in $\Sigma \cup\bsig$ is 
a $\cV_h$-regular embedding $\bar\cG \subset \Sigma \cup\bsig$. 
\end{definition}

\begin{figure}[h]
 \centering
     \begin{minipage}[t]{.9\textwidth}
      \centering
\includegraphics[angle=0, width=13cm, height=2.5cm]{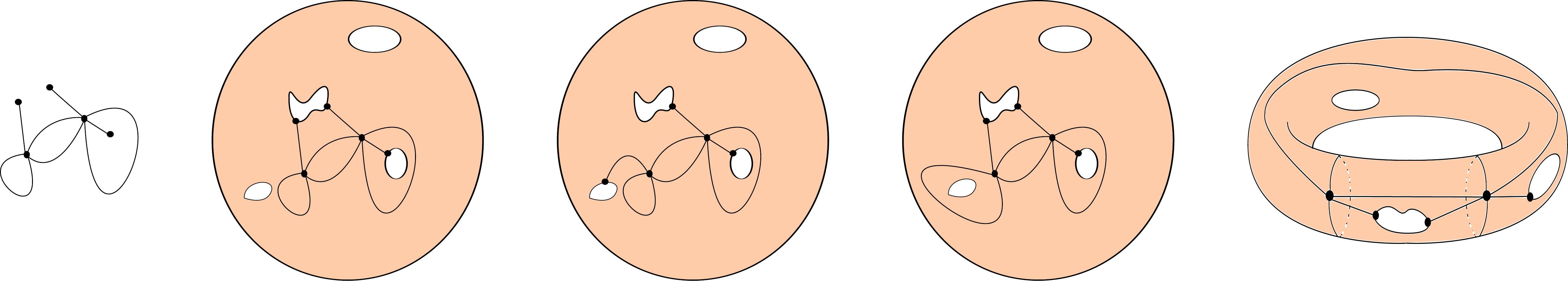}
\caption{ {\small The completed graph of the HEG of Figure \ref{fig:heg} and
some possible regular embeddings of it in a punctured sphere and a 
punctured torus. }}
\label{fig:regemb}
\end{minipage}
\end{figure}

Regular embeddings of (colored) graphs prove to be useful
in the context of graph encoded manifolds, see for
instance the work by Gagliardi in \cite{gagliardi}
and by Bandieri et al. in \cite{bandi}. 
Note that we choose to perform the cellulation on 
$\wsig$ to avoid subtleties induced by the natural topology
of $\Sigma \cup \bsig$.   
The last condition on the embedding, 
i.e.  $\bsig \cap \phi(\bar\cG)= \phi(\cV_h)$,
means that we require that the leaves in $\bar\cG$ obtained by completing the half-edges end on the boundary $\bsig$. 
Examples of regular embeddings for the completed graph
of Figure \ref{fig:hegstar} have been given 
in Figure \ref{fig:regemb}.

\begin{proposition}\label{prop:one}
Let $\cG_h$ be a HEG, $\bar\cG$ its  completed graph and $\cG$ its underlying graph. 
If  $\bar\cG\subset \Sigma \cup\bsig$ is a regular embedding, then there exists a cellular embedding $\cG\subset\wsig$.
\end{proposition}
\proof
A regular embedding  $\bar\cG\subset \Sigma \cup\bsig$ extends to a cellular embedding of $\bar\cG$ in $\wsig$  by extension of the codomain and keeping the cell decomposition of $\wsig-\phi(\bar\cG)$. The restriction $\phi|_{\cG}$ is a continuous map from $\cG$ to $\wsig$ and $\phi|_{\cG}: \cG\rightarrow \phi(\cG)$ is an homeomorphism as a restriction of the homeomorphism $\phi|_{\bar\cG}: \bar\cG\rightarrow \phi(\bar\cG)$. Furthermore $\wsig-\phi(\bar\cG)$ is homeomorphic to $\wsig-\phi(\cG)$ because $\bar\cG$ and $\cG$ have same cycles in $\wsig$.

\qed

\begin{definition}[HEG cellular embedding]\label{def:celheg}
A HEG $\cG_h$ is cellularly embedded in a punctured surface $\Sigma$ if and only if there is a  regular  embedding $\bar\cG\subset \Sigma \cup \bsig$. We denote 
the HEG cellular embedding in $\Sigma$ by $\cG_h \subset \Sigma$. 
\end{definition}
The following proposition holds. 
\begin{proposition}\label{prop:cgsigm}
Let  $\cG_h$ be a  HEG cellularly embedded in a punctured surface $\Sigma$.
Then its underlying graph  $\cG$ is cellularly embedded in $\Sigma$.  
\end{proposition}
\proof 
In the above notations, consider the regular embedding $\bar\cG\subset \Sigma \cup \bsig$ associated with $\cG_h \subset \Sigma$.
By Proposition \ref{prop:one}, $\cG$ is cellularly embedded in $\wsig$ and $\wsig -\phi(\cG)$  is disjoint union of 2-cells.
Let us denote $\wsig=\Sigma\cup(\cup_iD_i)$ the surface obtained by capping off $\Sigma$ where the $D_i$'s are spaces homeomorphic to closed discs. The $D_i$'s do not intersect $\phi(\cG)$ then $\wsig -\phi(\cG)=(\Sigma -\phi(\cG))\cup(\cup_iD_i)$ equals a  union of 2-cells.  Thus $\Sigma -\phi(\cG)$ is equal to a disjoint union of 2-cells possibly with punctures introduced by the $D_i$'s.

\qed

 \begin{definition}\label{def:equihegcell}
Let $\Sigma$ and $\Sigma'$ be two surfaces with punctures. 
Two cellularly embedded HEGs $\cG_h\subset \Sigma$ and $\cG'_{h'}\subset \Sigma'$ are equivalent, if their corresponding regular embeddings $\bar\cG \subset \Sigma \cup \bsig$ and $\bar\cG' \subset \Sigma' \cup \bsig'$ are equivalent, that is  
if there is a homeomorphism $\eta: \Sigma\cup \bsig \rightarrow \Sigma' \cup \bsig'$ with the property that  $\eta|_{\bar\cG}: \bar\cG \rightarrow \bar{\cG'}$ is an  homeomorphism.
\end{definition}

\begin{figure}[h]
 \centering
     \begin{minipage}[t]{1\textwidth}
      \centering
\includegraphics[angle=0, width=7cm, height=1.7cm]{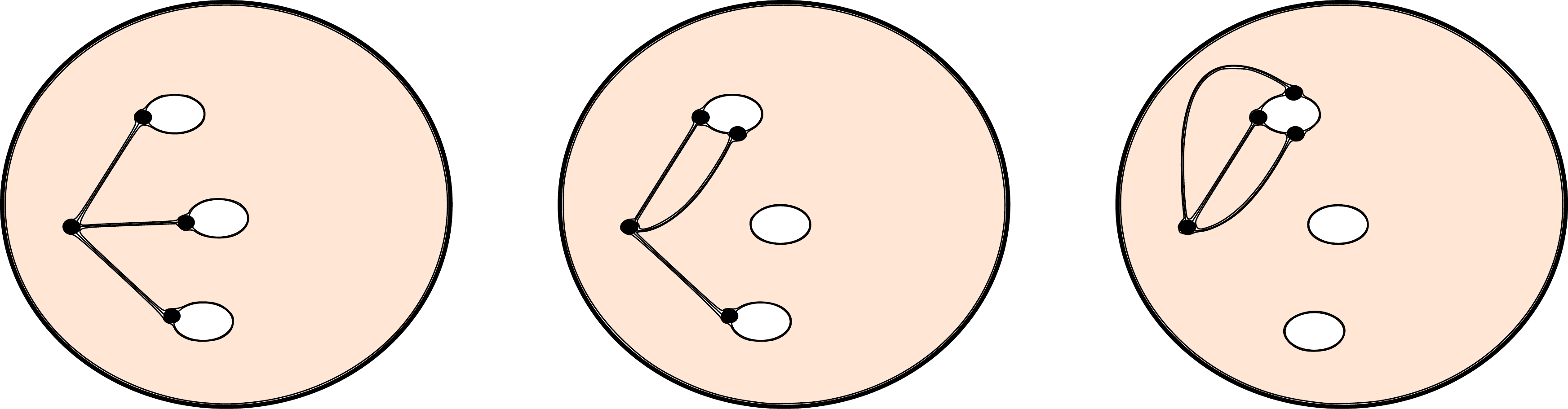}
\caption{ {\small Equivalent regular cellular embeddings
of the same completed graph of a HEG (a vertex with 3 half-edges) in the sphere with 3 punctures.}}
\label{fig:regembambi}
\end{minipage}
\end{figure}

A crux remark is that, while 
equivalence of cellular embeddings of graphs  
in (punctured) surfaces in the sense of Definition \ref{def:celemb} ensures that the number of 2-cells of the decomposition
is the same,  
 the equivalence of cellular 
embeddings of HEGs in punctured surfaces does not
anymore guarantees this property. See Figure \ref{fig:regembambi} for 
a simple illustration. 
This is source of ambiguities when we will 
seek equivalence between HEG cellular embeddings and HERGs
in the next paragraph. 
More restrictions on Definition \ref{def:equihegcell} could be discussed. 
For example, one could demand that the number of 2-cells should be the same after the cellulations of the two punctured surfaces $\Sigma$ and $\Sigma'$ to achieve equivalence of HEG cellular
embeddings. However, one can check this particular restriction will 
not lift the above mentioned ambiguity.

\medskip

\noindent{\bf Half-edge ribbon graphs (HERGs) -}
The class of ribbon graphs extends to the class of HERGs with the introduction of half-edges which are now ribbon-like. Half-edges of HERGs will be called half-ribbons (HRs). A HR is a ribbon incident to a unique vertex of a ribbon graph by a unique line segment on its boundary and without forming a loop. 

\begin{figure}[h]
 \centering
     \begin{minipage}[t]{1\textwidth}
      \centering
\includegraphics[angle=0, width=4cm, height=2cm]{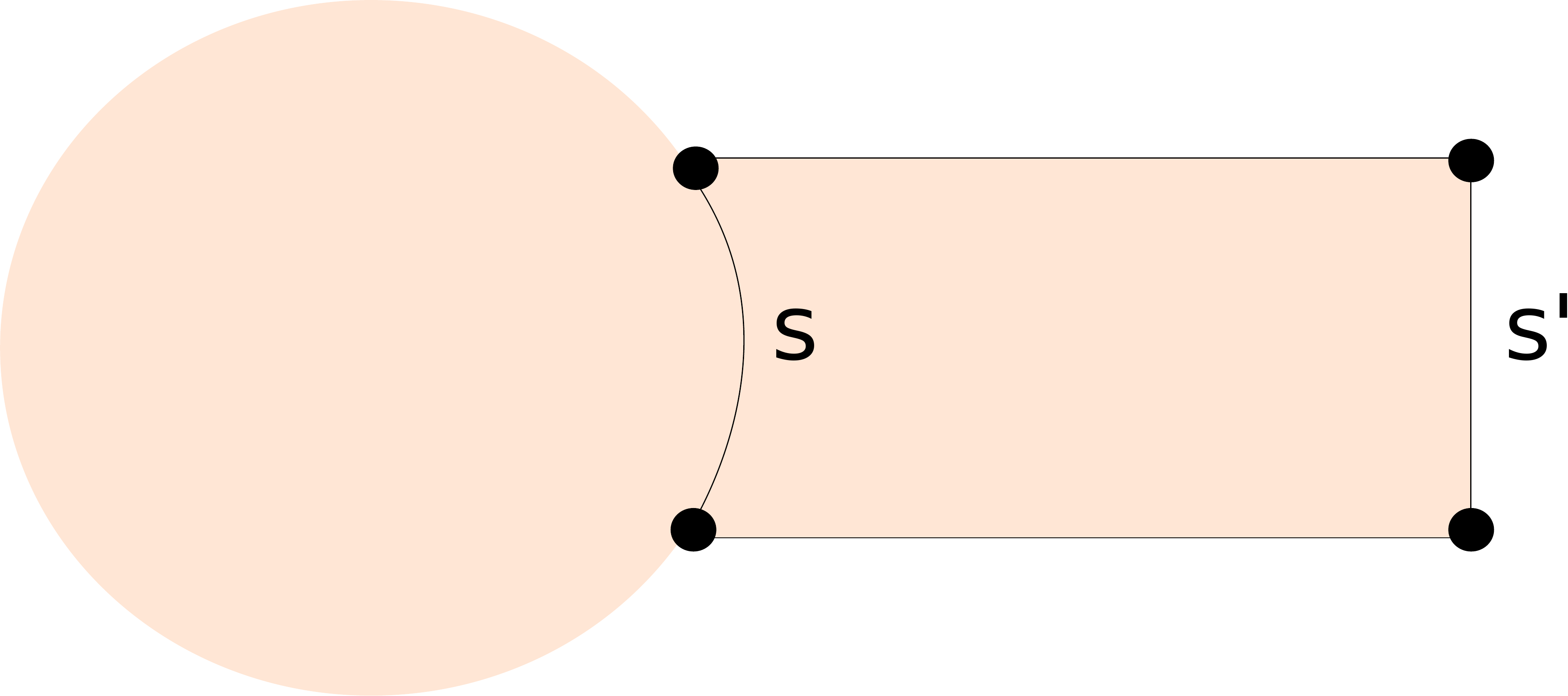}
\caption{ {\small A HR (with rectangular shape) incident to a vertex disc; the segment $s$ in contact with the vertex and $s'$ the external segment of the HR.}}
\label{fig:hr}
\end{minipage}
\end{figure}

If topologically, ribbons are discs, for combinatorial purposes, we regard a HR as a rectangle rather than a disc. To  achieve this, we introduce 4 distinct marked points at the boundary of a disc. A HR is incident to a vertex along a unique boundary arc $s$ lying between two 
successive of these marked points.
The segment parallel to $s$ is called \emph{external segment}. The end-points of any external segment are called \emph{external points} of the HR.  Figure \ref{fig:hr} illustrates a HR incident to vertex disc.

\begin{definition}[HERG]\label{def:herg}
A HERG  $G(V,E,H)$, or simply $G_{H}$, 
is a ribbon graph $G(V,E)$  with a set $H$ of HRs together with an incidence relation which associates each HR with a unique vertex. The ribbon graph $G$ is called the underlying ribbon graph of $G_H$. There is an
underlying HEG $\cG(V,E,H)$, denoted $\cG_H$, obtained from $G_{H}$ by keeping its
vertex, edge and half-edge sets and their incidence relation. 
\end{definition} 

Note that, as far as topology is concerned and as surfaces with boundary, HERGs are homeomorphic to ribbon graphs. HERGs have however richer combinatorial properties: we will use 
the modifications introduced by HRs to encode topological information such as punctures of surfaces in which cellular embeddings will be made.

\begin{figure}[h]
 \centering
     \begin{minipage}[t]{1\textwidth}
      \centering
\includegraphics[angle=0, width=5cm, height=1.7cm]{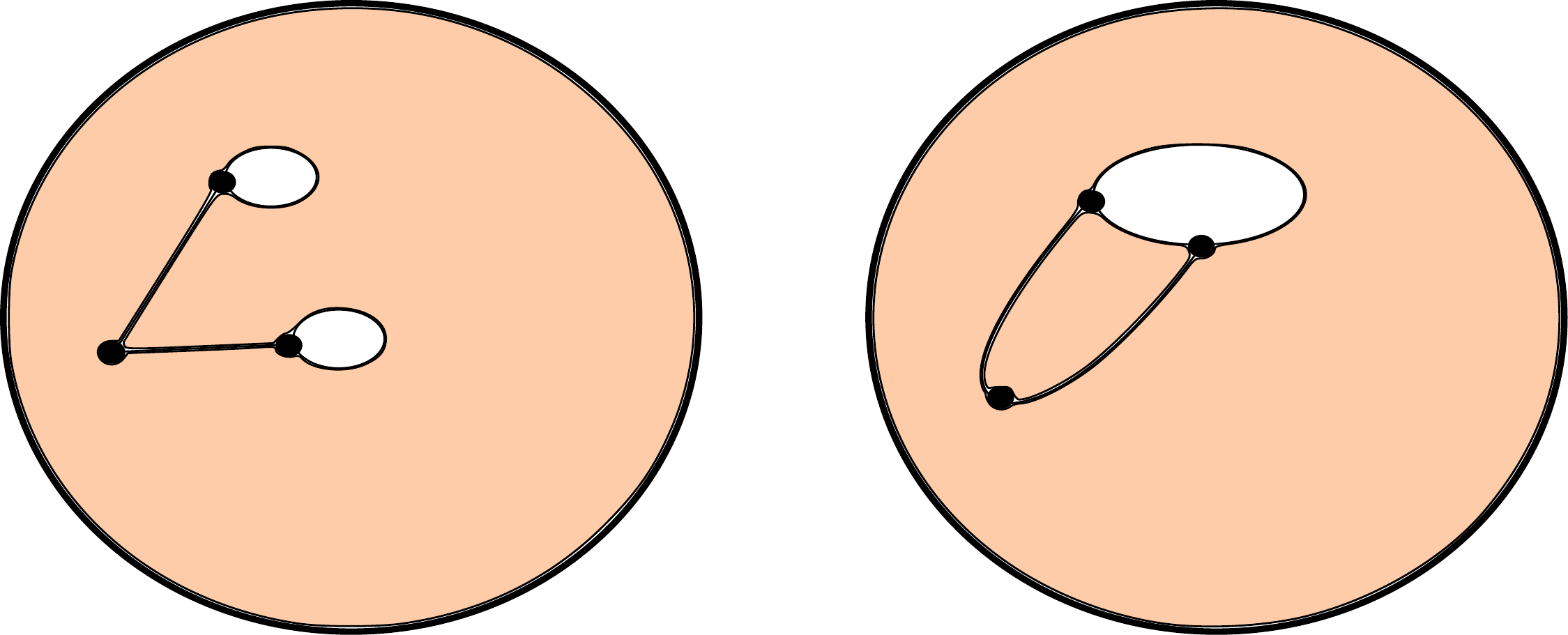}
\caption{ {\small Unequivalent cellular embeddings
of the same HEG (a vertex with 2 half-edges) in different punctured spheres giving rise to the same HERG.}}
\label{fig:regemuneq}
\end{minipage}
\end{figure}

Consider a HEG $\cG_h$ cellularly embedded in a punctured surface $\Sigma$. 
A HERG $G_H$ representation of the embedding of $\cG_h$ is obtained by 
taking  a small neighborhood band around $\cG_h$ in $\Sigma$. 
The set $H$ of HRs of $G_H$ is the set of neighborhoods of half-edges of $h$. 
Thus, to identify a HERG from a cellularly embedded HEG, the procedure is straightforward. 
However, the reverse procedure is ambiguous in several ways: 
there are several cellular embeddings of the same HEG in punctured surfaces (homeomorphic or not) which  all have the same HERG via the above procedure. Starting with a fixed (up to homeomorphism)
punctured surface, the cellular embeddings yielding the same HERG are those equivalent in the sense of Definition \ref{def:equihegcell}
(see, again, Figure \ref{fig:regembambi}). 
Without further assumptions, 
there is no criteria to lift the ambiguity, in other words, 
any representative in the class of equivalent cellular embeddings can be used to represent the HERG. 
The trouble becomes more apparent if we work 
with non homeomorphic surfaces: 
there are indeed  cases of nonequivalent cellular embeddings of the same HEG giving rise to the same HERG.
An example has been given in Figure \ref{fig:regemuneq}.  
To find a procedure which uniquely selects the cellular embedding $\cG_h \subset \Sigma$ from a HERG $G_H$, we need either restrictions on the definition of cellular embeddings of HEGs or some minimal requirements to choose one among those embeddings. Using the combinatorics of HERGs,
this problem has  at least two solutions. One of the prescriptions is in some sense
\emph{minimal}, the other maximal, and so could either be adopted by convention.

Given a HERG $G_H$, our goal is to construct a cellular embedding of some HEG in some punctured surface $\Sigma$ which obey unambiguously Definition \ref{def:celheg}.

It is not complicated to extend the completing procedure for HEG 
given by Definition \ref{completeprun} to HERGs. Consider a HERG $G_H$ and 
$V_H$ a set of discs such that $|V_H|=|H|$. 
We introduce $\bar G(V\cup V_H,E\cup H)$ the \emph{completed ribbon graph} 
obtained from $G_H$ by adding new vertices, elements of the set $V_H$,   
such that each HR in $H$ becomes a ribbon edge incident to a
unique vertex of $V_H$.

Since $\bar G$ is an ordinary ribbon graph,
the standard procedure to find the corresponding cellular embedding applies to it: we can 
find a graph $\cG_0$ which is cellularly embedded
in some closed connected compact surface $\Sigma_0$ of minimal genus. $\Sigma_0$ is
the capping off of $\bar G$ and
its genus is that of $\bar G$. 
A moment of thought, one easily 
realizes that $\cG_0$ is the completed graph $\bar\cG$ of $\cG_H$, the latter being 
the underlying HEG of $G_H$. 
The next move is to produce a cellular embedding of $\cG_H$ in some punctured surface $\Sigma$
obtained from $\Sigma_0$. This introduces two unknown data: the number $C_\bsig$ of boundary
circles in $\Sigma$ and the distribution of the half-edges of $H$ on these circles. 
As stated, this problem becomes purely combinatorial.

We now use the fact that, in HERGs, we distinguish several types
of boundary components \cite{avohou}.

\begin{figure}[h]
 \centering
     \begin{minipage}[t]{1\textwidth}
      \centering
\includegraphics[angle=0, width=4.5cm, height=2.8cm]{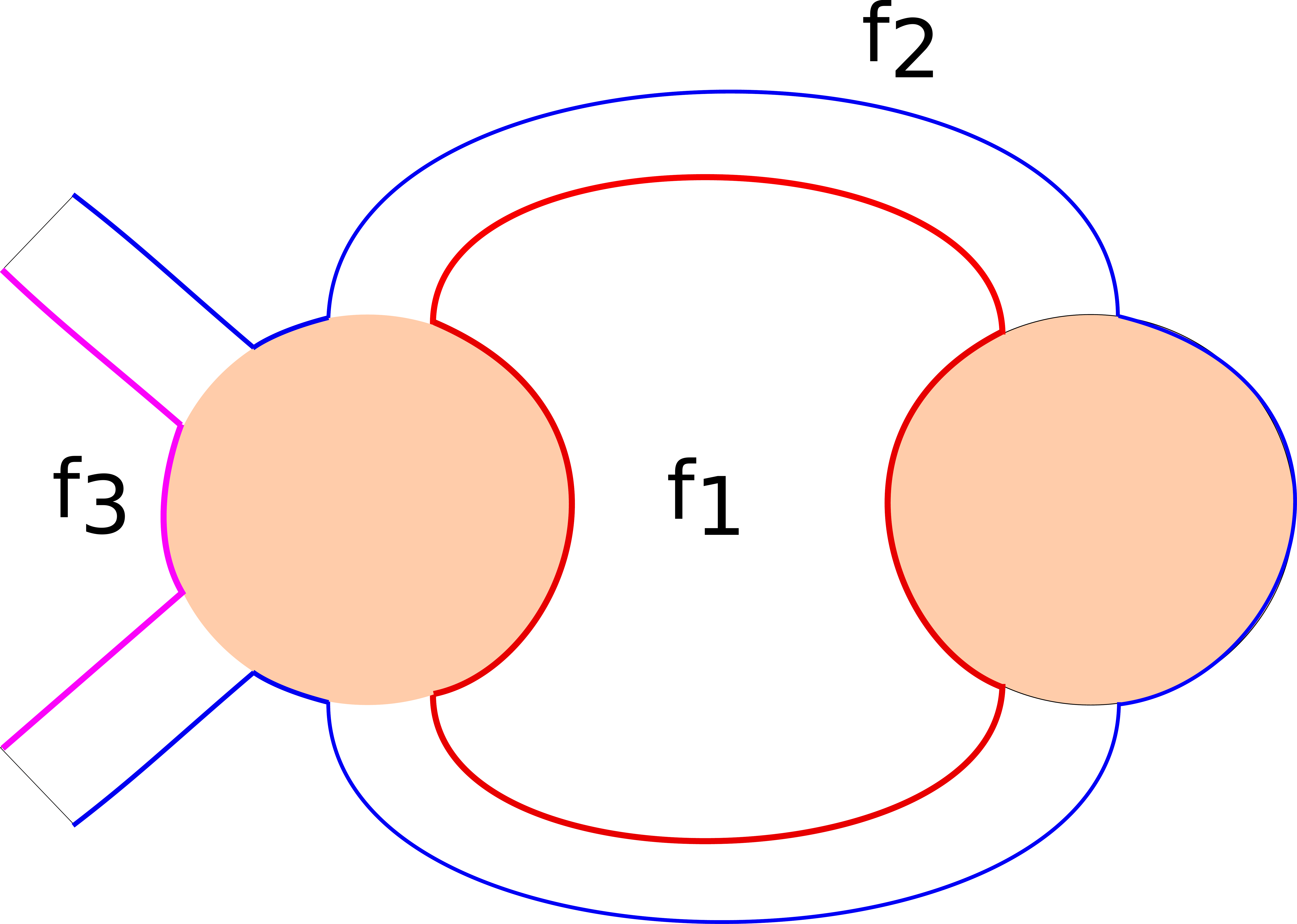}
\caption{ {\small A HERG with a closed face $f_1$ and
two open faces $f_2$ and $f_3$.}}
\label{fig:facinout}
\end{minipage}
\end{figure}

\begin{definition}[Closed and open faces]
\label{faces}
Consider  a HERG $G_{H}$.
A \emph{closed face} is a boundary  component of $G_H$  which never intersects any external segment of an HR.  The set of closed faces is denoted $\cF_{\inter}$.  An \emph{open  face}  is a boundary arc between an external point of some HR and another external point  without intersecting any external segment of an HR.  
 The set of open faces is denoted  $\cF_{\ext}$. 
 The set of faces $\cF$ of $G_H$ is defined by 
$\cF_{\inter} \cup \cF_{\ext}$. (See illustrations on Figure \ref{fig:facinout}.)
\end{definition}

We complete Definition \ref{faces}  by identifying
a new type of boundary component: 
 \begin{definition}[External cycles]
A boundary component of $G_H$ obtained by following alternatively external faces and external segments of the HRs is called {\it external cycle}.
\end{definition}
For a HERG $G_H$, following external cycles, we obviously have 
\bea\label{fach}
|\cF_{\ext}| = |H|. 
\eea

 External cycles form connected components of a 2-regular graph called in the boundary graph of $G_H$ in  \cite{Gurau:2009tz}.
 
The Euler  characteristic $\chi(\bar G)$ of the completed ribbon graph $\bar G$ of $G_H$ is,  using similar notations as above,  
\bea\label{eq:eulerheg}
\chi(\bar G) &=& v(\bar G) - e(\bar G) + f (\bar G)\cr
&=&   v(G_H)  - e(G_H)  + f_{\inter}(G_H)+C_\ext(G_H)
\eea
where $f(\bar G)$, the number of faces of $\bar G$, equals the number of closed faces $f_{\inter}(G_H)= |\cF_{\inter}(G_H)|$ plus the number of external cycles $C_\ext(G_H)$  of $G_H$. We note
also that $\chi(\bar G)=\chi(G)$, where
$G$ is the underlying ribbon graph of $G_H$.
Indeed, these ribbon graphs have the same number of faces, since a face in $\bar G$ corresponding to $C_\ext(G_H)$ deforms uniquely onto a face of $G$. Thus,  equivalently, 
the underlying graph $\cG$ of $G$ can be used
 to define the same surface $\Sigma_0$,
because $\cG \subset \Sigma_0$ is a cellular embedding 
equivalent to $G$.  We will use this remark in the following when
we will distinguish different cases of
HEGs cellular embeddings. 

We realize that a cellular embedding
in the sense of Definition \ref{def:celheg},
or equivalently a regular embedding
$\bar\cG \subset \Sigma \cup \bsig$,
gives us a constraint on $C_\bsig$ the number  of boundary components of the surface $\Sigma$ that
we are seeking. Mapping half-edges on 
the boundary circles, we  infer that  
$C_\bsig\geq C_\ext(G_H)$. Indeed, 
we recall that the set of half-edges is in one-to-one
correspondence with $H$ the set of HRs of $G_H$. 
$H$ is partitioned in $C_\ext(G_H)$ parts. 
On the other hand, the leaves in $\bar\cG$ 
corresponding to half-edges (and so to HRs) intersect necessarily a boundary circle. The inequality
therefore holds.  
The punctured surface is  obtained
$\Sigma$ after removing $C_\bsig\geq C_\ext(G_H)$ boundary circles  in the surface $\Sigma_0$. 
Finally, $\cG_h$ can be cellularly embedded in 
 $\Sigma$ 
or in any other punctured surface $\Sigma'$ with same genus and $C_\bsig'\geq C_\ext(G_H)$.

We make another observation: 
 if we request that any boundary circle on the surface must be intersecting
a vertex of the completed graph $\bar \cG$ for 
any regular cellular embedding $ \bar \cG \subset \Sigma \cup \bsig$, then we have also an upper bound on the 
number of boundary circles in the surface such that 
\bea
C_\ext(G_H) \leq C_\bsig \leq |H| \,. 
\eea
We now discuss two particular prescriptions
specializing the HEG cellular embedding
$\cG_H \subset \Sigma$ and their consequences.

$\bullet$ Assume $C_\bsig = C_\ext(G_H)$ which is the minimum number of boundary circles of $\Sigma$ to
ensure that there is a cellular embedding $\cG_H\subset\Sigma$ corresponding to the HERG $G_H$.
  Consider $\cG$ the underlying
graph of $\cG_H$. Proposition \ref{prop:cgsigm} instructs us that there is a cellular
embedding  $\cG \subset \Sigma$. 
Observe that $\Sigma-\phi(\cG)$ splits in two sets:  the set $O$ of 2-cells and the remaining set $O'$
of 2-cells with punctures. 
As explained previously, $G$ has a genus determined by \eqref{eq:eulerheg}, and therefore  the number of connected components of $\Sigma-\phi(\cG)$ is $f_{\inter}(G_H) + C_\ext(G_H)$. 
There are one-to-one correspondences, 
on one side, between the set of closed faces of $G_H$ and $O$ and, on 
the other side, between the set of 
external cycles of $G_H$ and $O'$. 
The equality $C_\bsig = C_\ext(G_H)$ simply reveals 
that each element of $O'$ has a single puncture.   
Furthermore, consider the cellulation $\Sigma-\phi(\cG_H)$ which gives a disjoint union of 2-cells. The number of such topological discs is equal to $f_{\inter}(G_H)+f_{\ext}(G_{H})$ where $f_{\ext}(G_H)$ is the number of external faces of $G_{H}$.
Indeed, it must be obvious from 
the previous comments that 2-cells of $O$ which are
again in $\Sigma-\phi(\cG_H)$ coincide
with closed faces of $G_H$ and their number corresponds to 
$f_{\inter}(G_H)$. 
Consider now the remaining set $O'$ of the cellulation $\Sigma - \phi(\cG)$. 
Each 2-cell of $O'$ has a single puncture and 
so a single boundary circle in its interior.
That 2-cell with one puncture may split after removing half-edges of $\cG_H$ in the cellulation $\Sigma-\phi(\cG_H)$. The number of parts of this splitting is precisely the number of half-edges
of its corresponding external cycle. 
Finally, a half-edge uniquely corresponds  to a HR and, 
in each external cycle of $G_H$, the number of open
faces is equal to the number of HR, see \eqref{fach}.  

There is an alternative prescription leading
to another unambiguous construction of the HEG cellular embedding
for a given HERG.

$\bullet$   Suppose that for each half-edge we assign a puncture
in the surface, which means that
we construct a surface such that $C_\bsig= |H|$. This leads also
to an unambiguous situation. Indeed, after capping off $\bar G$, the completed
graph of $G_H$, we are first led
to a closed surface $\Sigma_0$. Then, we prune $\bar G$
with respect to $V_H$ in $\Sigma_0$ (or simply remove the leaves of $\cG_H$ corresponding to $H$). 
Nevertheless, compared to the previous prescription, the number of boundary circles can not be minimal.

We finally introduce the notion of minimal/maximal HEG  cellular embedding corresponding to HERGs. 

\begin{figure}[h]
 \centering
     \begin{minipage}[t]{1\textwidth}
      \centering
\includegraphics[angle=0, width=8cm, height=3cm]{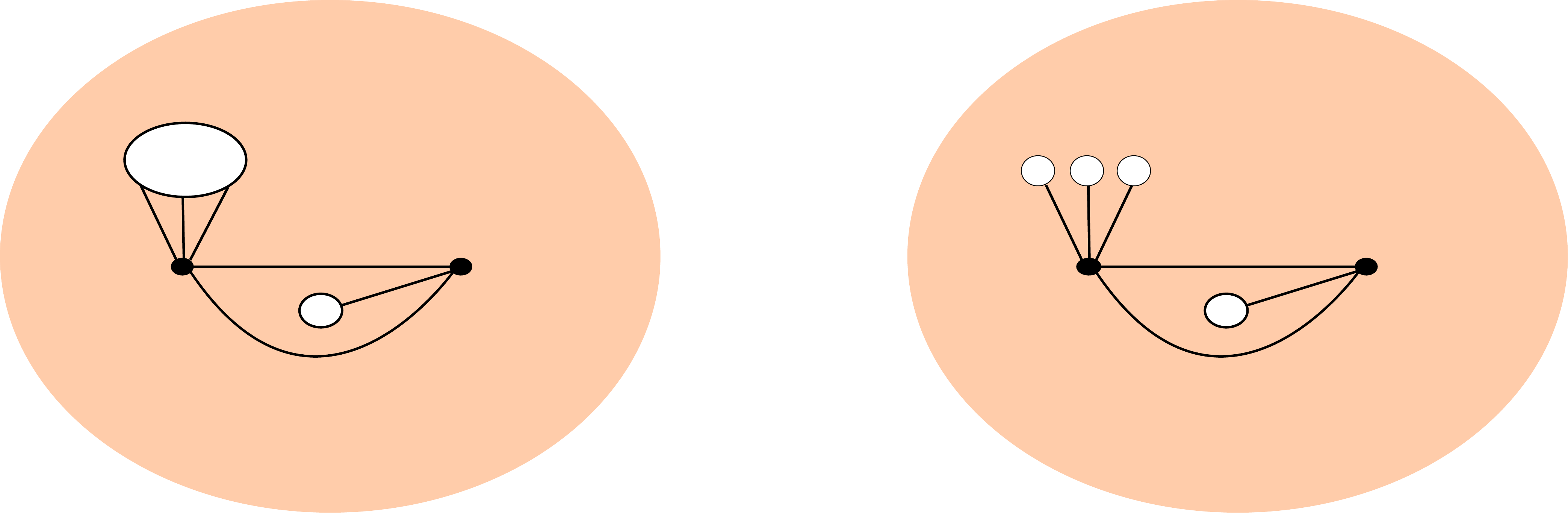}
\caption{ {\small Proper (left) and $h$-proper HEG cellular embeddings
on the punctured sphere. 
 }}
\label{fig:properh}
\end{minipage}
\end{figure}

\begin{definition} 
Consider $\cG_h$ a HEG and $\Sigma$ a punctured surface. 
A HEG cellular embedding $\cG_h\subset \Sigma$ is called proper if and only if the number of external cycles
of the HERG $G_H$ generated by the embedding 
equals the number of punctures of $\Sigma$. 

A HEG cellular embedding $\cG_h\subset \Sigma$ is called $h-$proper if and only if the number of half-edges equals the number of punctures of $\Sigma$
 and any boundary circle of $\Sigma$ intersects $\bar\cG$. 
\end{definition}
Figure \ref{fig:properh} gives some examples of proper and $h$-proper HEG 
cellular embedding on the punctured sphere.  

The following statement therefore holds: 
\begin{theorem}\label{theo}
A HERG corresponds to unique  proper ($h$-proper)
HEG  cellular embedding in a punctured surface with minimal genus and minimal (maximal) number of punctures.   

\end{theorem}

\begin{figure}[h]
 \centering
     \begin{minipage}[t]{1\textwidth}
      \centering
\includegraphics[angle=0, width=10cm, height=5cm]{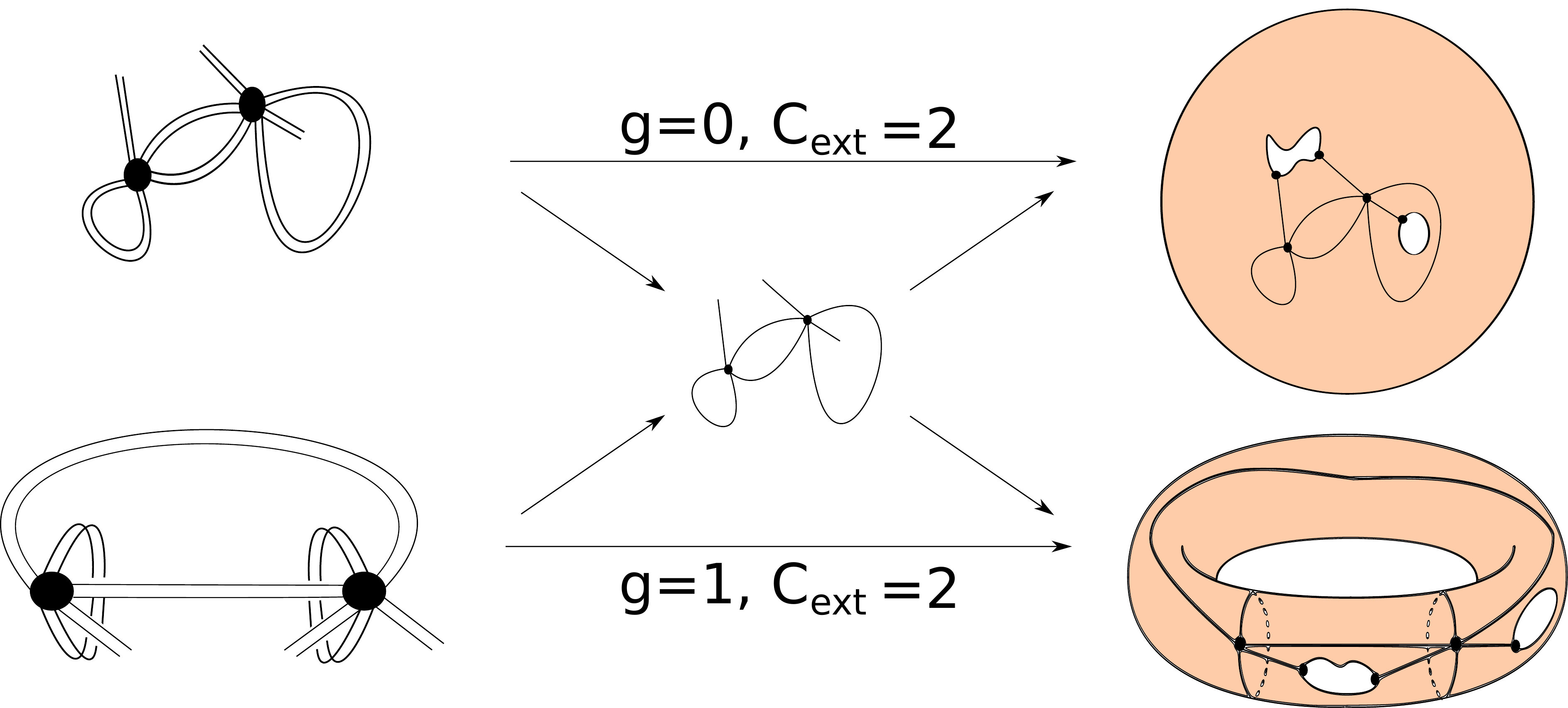}
\caption{ {\small Two HERGs and their HEG proper cellular embedding. 
 }}
\label{fig:regeminemb}
\end{minipage}
\end{figure}

Figure \ref{fig:regeminemb} shows HERGs with their unique
(up to homeomorphism) HEG proper cellular 
embedding. For $h$-proper HEG cellular embeddings corresponding
to the same HERG, one must put boundary discs for 
each half-edge in the surface. Finally, to close this section, it is obvious that the above construct reduces to usual graph cellular embedding on closed surfaces when the HEG does not have any half-edges and so no punctures are needed on the surface.  

\section{Geometric duality for HEG cellular embeddings}
\label{sect:geomdual}

In this section, we generalize the geometric duality of graphs cellularly embedded in surfaces to cellular embedded HEGs in punctured surfaces. Although our illustrations are only made on the 2-sphere, the duality is valid on any punctured surface.

\medskip

\noindent{\bf Geometric duality for graph cellular embeddings -}
Let $\cG(\cV,\cE)$ be a graph and $\cG\subset\Sigma$ be a cellular embedding of $\cG$ in a closed connected compact surface $\Sigma$. The geometric dual  $\cG^*\subset\Sigma$ of  $\cG\subset\Sigma$ 
is the cellular embedding  in $\Sigma$ of the graph $\cG^*(\cV^*,\cE^*)$ obtained by inserting one vertex in each of the faces of $\cG$ and embedding an edge of $\cE^*$ between two of these vertices if the faces of $\cG$ where they belong are adjacent. Then an edge of $\cE^*$ crosses the corresponding edge of $\cE$ transversely. An edge of $\cE^*$ forms a loop if it crosses an edge of $\cE$ incident to only one face of $\cG$. Hence the set of vertices of $\cG^*$ are  in one-to-one correspondence with the set of faces of $\cG$,
and $|\cE^*|=|\cE|$, $|\cV^*|=f(\cG)$ where $f(\cG)$ is the number of faces of $\cG$. Thus, we have $\gamma(\cG)=\gamma(\cG^*)$.

\medskip 

\noindent{\bf Geometric duality for graph cellular embeddings  in punctured surfaces -}
If $\cG\subset\Sigma$ is a cellular embedding in a punctured surface $\Sigma$, then we can also construct $\cG^*\subset\Sigma$ a cellular embedding in $\Sigma$ by the same recipe developed in the previous paragraph. The construction of $\cG^*$ simply  
avoids the punctures and lies in the interior of the punctured surface. 
Once again, we can cap off the surface,
determine the geometric dual and insert back the punctures on the surface.

\medskip

\noindent{\bf Geometric duality for HEG cellular embeddings in a punctured surfaces -}
Consider a HEG cellular embedding $\cG_h\subset \Sigma$ and
 $\bar\cG\subset \Sigma\cup \bsig$ its associated regular embedding
and $\cG \subset \Sigma$ the cellular embedding of its underlying
graph (Proposition \ref{prop:cgsigm}). 
We want to define the dual $(\cG_h)^*\subset \Sigma$
of $\cG_h\subset \Sigma$. 
 The first track is to construct the geometric dual $(\bar\cG)^* \subset\Sigma\cup \bsig$ of
the  regular embedding $\bar\cG \subset\Sigma\cup \bsig$
of the completed graph $\bar \cG$
and, then, operate on  $(\bar\cG)^*$ to identify what 
 $(\cG_h)^*$ could be. 

We henceforth work under two conditions: 
(1) after the embedding, all boundary circles of the surface intersect
at least one vertex of the completed graph; 
(2)  the dual of a 
HEG has the same number of half-edges of the HEG. The property (2) was indeed shown true in \cite{krf} in the case of duals HERGs. We would like to preserve this feature for
duals of cellularly embedded HEGs. 

$\bar\cG$ is a graph regularly embedded in a surface $\Sigma\cup\bsig$ 
and intersect the boundary $\bsig$. 
We cannot construct its dual 
$(\bar\cG)^*$ according to the previous paragraph by simply avoiding
the boundary $\bsig$. 
Constructing that dual, consider rather the 2-cells
intersecting boundary circles obtained from  the cellulation 
$\Sigma - \phi(\cG)$. Each of these 2-cells
might further split after removing the edges of $\bar\cG$
which are the completed of the half-edges of $\cG_h$. 
We call these 2-cells with punctures {\it external cycles}
of the HEG $\cG_h$.  For 2-cells without punctures,
constructing the dual graph remains the same
as in usual situation: dual vertices are defined by 1 vertex per such 2-cell and 
dual edges transversal to the edges of the 2-cell. 
We will therefore focus on the duality at the level of external cycles. 

\begin{figure}[h]
 \centering
     \begin{minipage}[t]{1\textwidth}
      \centering
\includegraphics[angle=0, width=7.5cm, height=3cm]{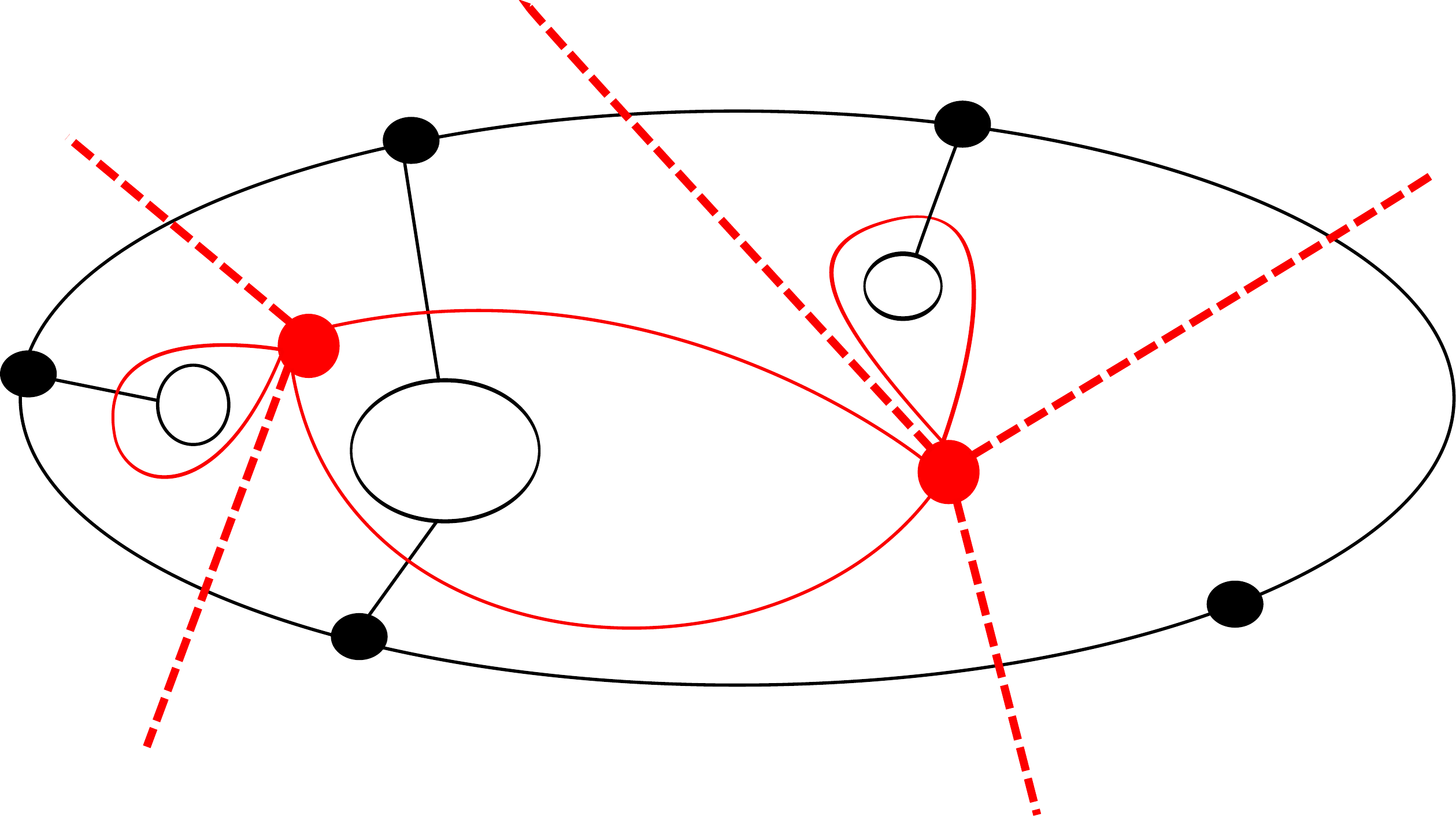}
\caption{ {\small 
Constructing the dual $\bar\cG^*$ (in red) at an external 
cycle with 5 vertices: with each boundary circle within the external cycle, there is a cycle 
associated in $\bar\cG^*$.  Dual edges are put in dash and are incident to other vertices.}}
\label{fig:dualext}
\end{minipage}
\end{figure}

Given condition (1), 
with  each boundary circle there is an associated cycle in $\bar\cG^*$
 (see Figure \ref{fig:dualext}). Either those cycles are loops and 
they become one-to-one with a subset of boundary circles or the cycles
are of length larger than 2.   

 We  modify $\bar\cG^*$ to define the dual of $\cG_h \subset \Sigma$. 
 To satisfy condition (2), working on a  $h$-proper HEG cellular embedding, 
  we can use a mapping: to each loop of the dual $\bar\cG^*$ lying in an external cycle of $\cG_h$, we associate a half-edge. The procedure becomes unambiguous because
all cycles in an external cycle of $\cG_h$ are loops. 
In other situations, we need more work. 

\begin{definition}[$h$-weak HEG cellular embedding]
Consider a HEG cellular embedding $\cG_h\subset \Sigma$. 
Then $\cG_h\subset \Sigma$  is called a $h$-weak  HEG cellular embedding if the edges
completing the half-edges of $h$ in  $\bar\cG\subset \Sigma\cup \bsig$, which are incident to the same boundary circle are also incident to the same vertex in $\bar\cG$ (or in $\cG_h$). 
\end{definition}

\begin{figure}[h]
 \centering
     \begin{minipage}[t]{1\textwidth}
      \centering
\includegraphics[angle=0, width=4cm, height=3cm]{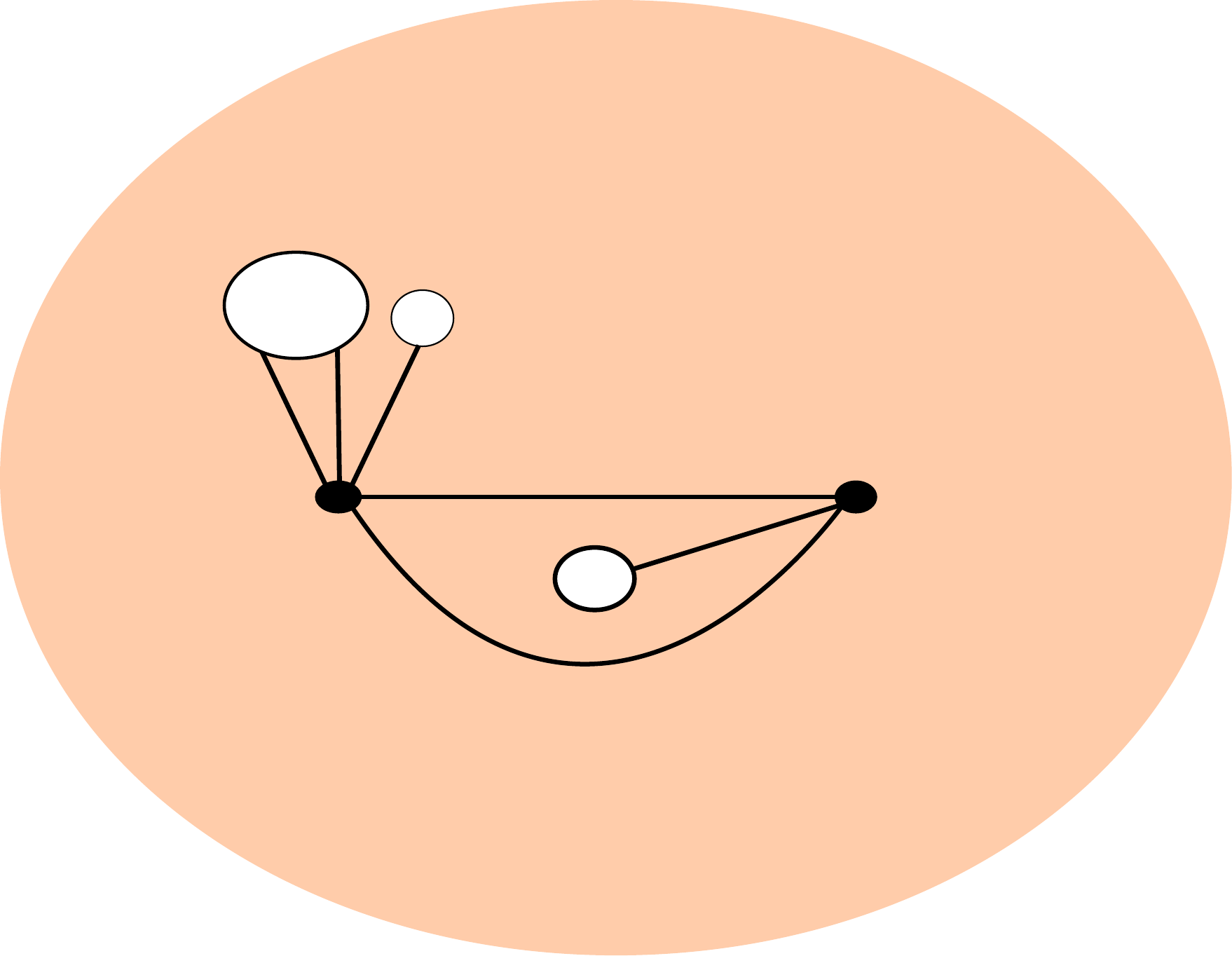}
\caption{ {\small 
A $h$-weak HEG cellular embedding on the punctured sphere.}}
\label{fig:hweak}
\end{minipage}
\end{figure}

Note that a $h$-proper HEG cellular embedding is  $h$-weak.  
A $h$-weak HEG cellular embedding has been illustrated in Figure \ref{fig:hweak}. In each external 
cycle, we note that there is a \emph{special vertex} of $\bar\cV^*$ of  $\bar\cG^*$,
which is of degree twice number of cycles  plus the number
of edges forming the external cycle. 
This vertex will be useful in the following operations. 

Let $\cG_h=\cG(\cV,\cE,h)$ be a HEG with 
completed graph $\bar\cG(\bar\cV=\cV\cup \cV_h, \bar\cE=\cE\cup h)$. 
Let $\cG_h \subset \Sigma$ be a $h$-weak HEG cellular embedding,
$\bar\cG\subset \Sigma\cup\bsig$ be its regular embedding. 
Consider $\bar\cG^*\subset \Sigma\cup\bsig$ the geometric dual of  $\bar\cG\subset \Sigma\cup\bsig$,  that we denote as a graph as  
 $\bar\cG^*(\bar\cV^*,\bar\cE^*)$. 
Let us denote $\cV_{\ext}^*$ the subset of vertices of $\bar\cV^*$ associated with 
a given external cycle of $\cG_h$ in $\Sigma$ and call $v^*$ the special vertex
of the external cycle. 
Consider $L_h$ the set of cycles formed with the dual of the edges belonging to 
$h \subset \bar\cE$.  
Each $l\in L_h$ is encircling
a boundary circle $c_l$ where ends corresponding edges $\{e_{l,k}\}$ elements of $h$. 
If the cycle $l$ is a loop, there is a single edge $e_l\in h$. All cycles $l$
are incident to $v^*$. The length of a cycle $l$ is denoted $|l|$. 
We want to regard $L_h$ as a subset of edges hence we write
$L_h \subset \bar\cE^*$; the set of vertices which forms $L_h$ 
%except $v^*$ 
is denoted $\cV(L_h)$. We have $|\cV(L_h)|+1=|\{e_{l,k}\}|=|\{v_{l,k}\}|$.

\begin{definition}[Grafting and grafted graph]
The \emph{grafting} operation on a cycle $l \in L_h$
circumventing a boundary circle $c_l$ and incident to $v^*$, consists 
in  removing all edges of $l$ and the vertices where these edges are incident except $v^*$, 
then inserting  $|l|$ embedded edges $e^*_{l,k}$, $k=1,\dots,|l|$, 
respecting the cyclic ordering around $v^*$ and keeping all remaining edges and vertices untouched.
The edges $e^*_{l,k}$ are incident to $v^*$ and to new vertices $v_{l,k}$ on $c_l$.
 
The \emph{grafted graph} $\cG_{1;L_h}=\cG_1(\cV_1,\cE_1)$ with respect to $L_h$ 
is the graph with vertex set 
 $\cV_1=\bar\cV^* \setminus \cV(L_h) \cup \{v_{l,k} \}_{l \in L_h}$ and edge 
$\cE_1 = \bar\cE^* \setminus L_h \cup \{e^*_{l,k}\}_{l \in L_h}$
obtained after  performing a sequence of grafting operations
 on $\bar\cG^*$,  for all   $l\in L_h$. 
\end{definition}

\begin{figure}[h]
 \centering
     \begin{minipage}[t]{1\textwidth}
      \centering
\includegraphics[angle=0, width=9cm, height=3.5cm]{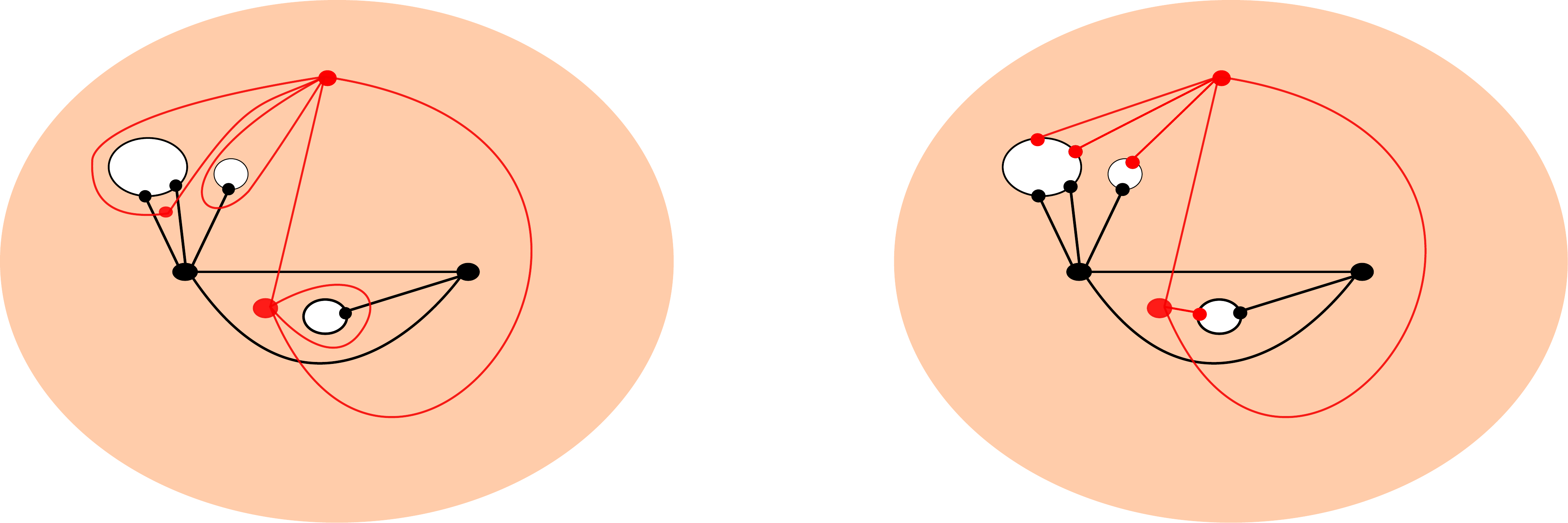}
\caption{ {\small 
The grafting of  graph cellular embedding (in red) on the punctured sphere.}}
\label{fig:hweakdual}
\end{minipage}
\end{figure}

The grafting operation on a graph cellular embedding is shown in Figure \ref{fig:hweakdual}. 
To simplify notations, we write $\{v_{l,k}\}_{l\in L_h}$ as $\{v_l\}$
and $\{e^*_{l,k}\}_{l \in L_h}$ as $\{e^*_l\}$.

\begin{theorem}\label{Propo:dual}
Let $\cG_h = \cG(\cV,\cE,h)$ be a HEG,
$\cG_h \subset \Sigma$ be a $h$-weak HEG cellular embedding,
$\bar\cG\subset \Sigma\cup\bsig$ be its regular embedding with 
$\bar\cG(\bar\cV=\cV \cup \cV_h,\bar\cE=\cE\cup h )$,
 $\bar\cG^*\subset \Sigma\cup\bsig$ be the geometric dual of $\bar\cG\subset \Sigma\cup\bsig$, with $\bar\cG^*(\bar\cV^*,\bar\cE^*)$,  and $L_h\subset \bar\cE^*$ be the set of
 cycles formed with edges duals to edges belonging to $h\subset \bar\cE$.

The grafted graph $\cG_{1;L_h}=\cG_1(\bar\cV^* \setminus \cV(L_h) \cup \{v_l \},\bar\cE^* \setminus  L_h \cup \{e^*_l\})$  with respect to $L_h$ defines a $\{v_l \}$-regular embedding $\cG_1\subset \Sigma \cup \bsig$. Furthermore,
pruning $\cG_1$ with respect to $ \{v_l\}$ 
defines a $h'(\{v_l\})$-weak HEG cellular embedding in $\Sigma$, where $h'(\{v_l\})$ is the set of half-edges resulting from the pruning of the edges of $\{e^*_l\}$. 
\end{theorem}  
\proof 
The fact that $\cG_1 \subset \Sigma \cup \bsig$ is a $\{v_l \}$-regular embedding 
can be easily shown: $\cG_h \subset \Sigma$ being $h$-weak,
then a subset $S$ of  leaves  (which correspond to a subset of completed half-edges in $\bar\cG$) ending on a boundary circle in $\bsig$ maps in the geometric dual $\bar\cG^*\subset \Sigma\cup\bsig$ to a unique cycle $l_S$ encircling that boundary circle. 
Note that for all $S$, all boundary circles are encircled.   
The grafting of these cycle  $l_S$, for all $S$, makes them a subset of  edges $\{e_{l}^*\}$ with end vertices $\{v_{l}\}$ intersecting all boundary circles. 

We concentrate on the second
statement. Consider the HEG obtained after  pruning $\cG_{1;L_h}$
that we denote $\und{\cG_{1;L_h}}_{\,h'(\{v_l\})}$. 
The vertex set of $\und{\cG_{1;L_h}}_{\,h'(\{v_l\})}$ is given by
$\und{\cV}=\bar\cV^*\setminus \cV(L_h)$, its edge set by 
 $\und{\cE} =\bar\cE^* \setminus  L_h$
and its half-edge set by $h'(\{v_l\})=\{e^*_l\}$. 
The incidence relation
between $\und{\cE}$ and $\und{\cV}$ and between $h'(\{v_l\})$ 
and $\und{\cV}$ can be easily inferred since they  are
 inherited from extension and restriction of the incidence relations 
 in $\cG_{1;L_h}$. 
Clearly, completing $\und{\cG_{1;L_h}}_{ \,h'(\{v_l\})}$ 
gives back $\cG_{1;L_h}$
and then we can call $\cG_1 \subset \Sigma\cup \bsig$
a regular embedding. Thus there is a HEG cellular embedding 
of the HEG $\und{\cG_{1;L_h}}_{ \,h'(\{v_l\})}$ in $\Sigma$. 
The property that this  HEG cellular embedding is $h'(\{v_l\})$-weak
follows again by construction:  the set of half-edges in $\und{\cG_{1;L_h}}_{ \,h'(\{v_l\})}$
is in one-to-one correspondence with $h$ in $\cG_h$,
 $|h'(\{v_l\})|=|\{v_l \}|=|h|$, there is conservation of half-edges
after the procedure. We then conclude to the result since, per external cycle where the boundary circles are, all half-edges are incident to the same special vertex. 

\qed

\begin{definition}[Geometric dual of a $h$-weak HEG cellular embedding]
Let $\cG_h \subset \Sigma$ be a $h$-weak HEG cellular embedding.
The geometric dual of $\cG_h \subset \Sigma$, denoted by $\cG^*_h\subset \Sigma$,
is  the $h$-weak HEG cellular embedding constructed previously
 $\und{\cG_{1;L_h}}_{ \,h'(\{v_l\})}\subset \Sigma$. 
\end{definition}

\begin{figure}[h]
 \centering
     \begin{minipage}[t]{1\textwidth}
      \centering
\includegraphics[angle=0, width=7.5cm, height=3cm]{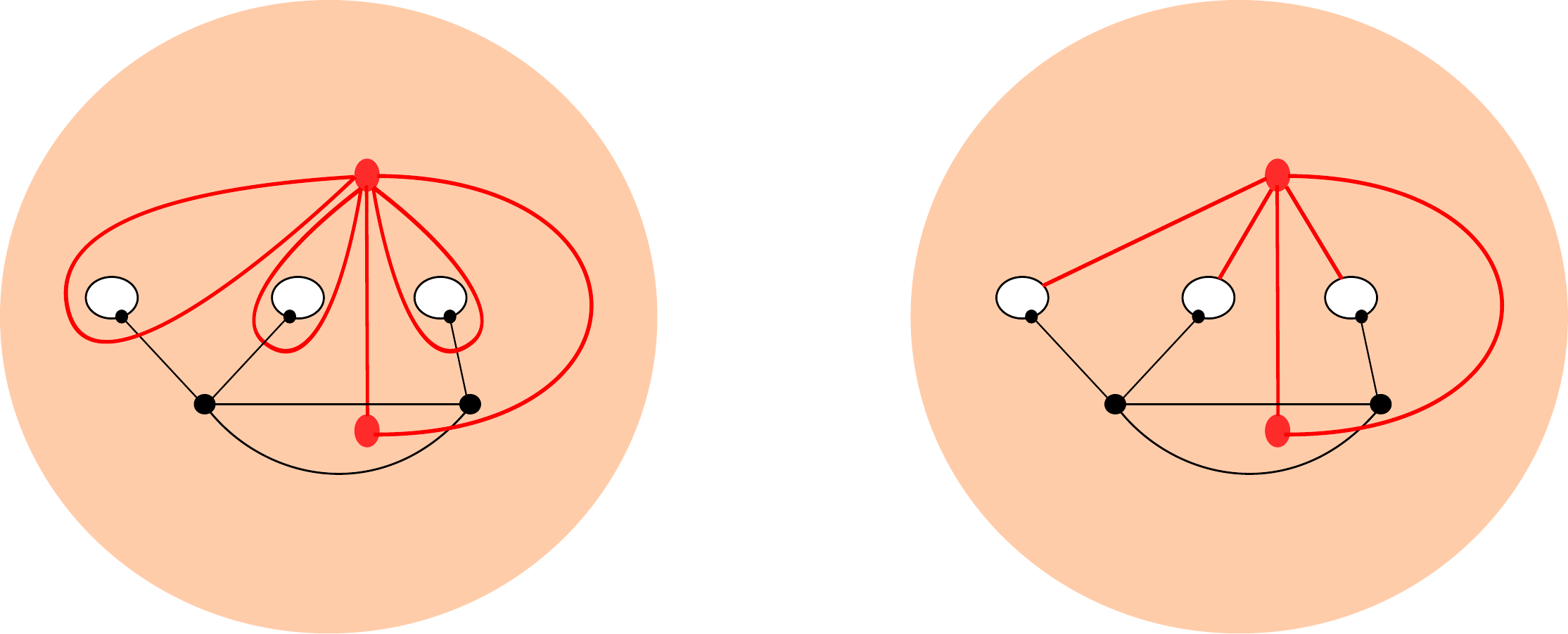}
\caption{ {\small The construction of the geometric dual of 
a $h$-proper HEG cellular embedding (in black) in a punctured sphere.}}
\label{fig:dualheg}
\end{minipage}
\end{figure}

\begin{figure}[h]
 \centering
     \begin{minipage}[t]{1\textwidth}
      \centering
\includegraphics[angle=0, width=4.5cm, height=3.5cm]{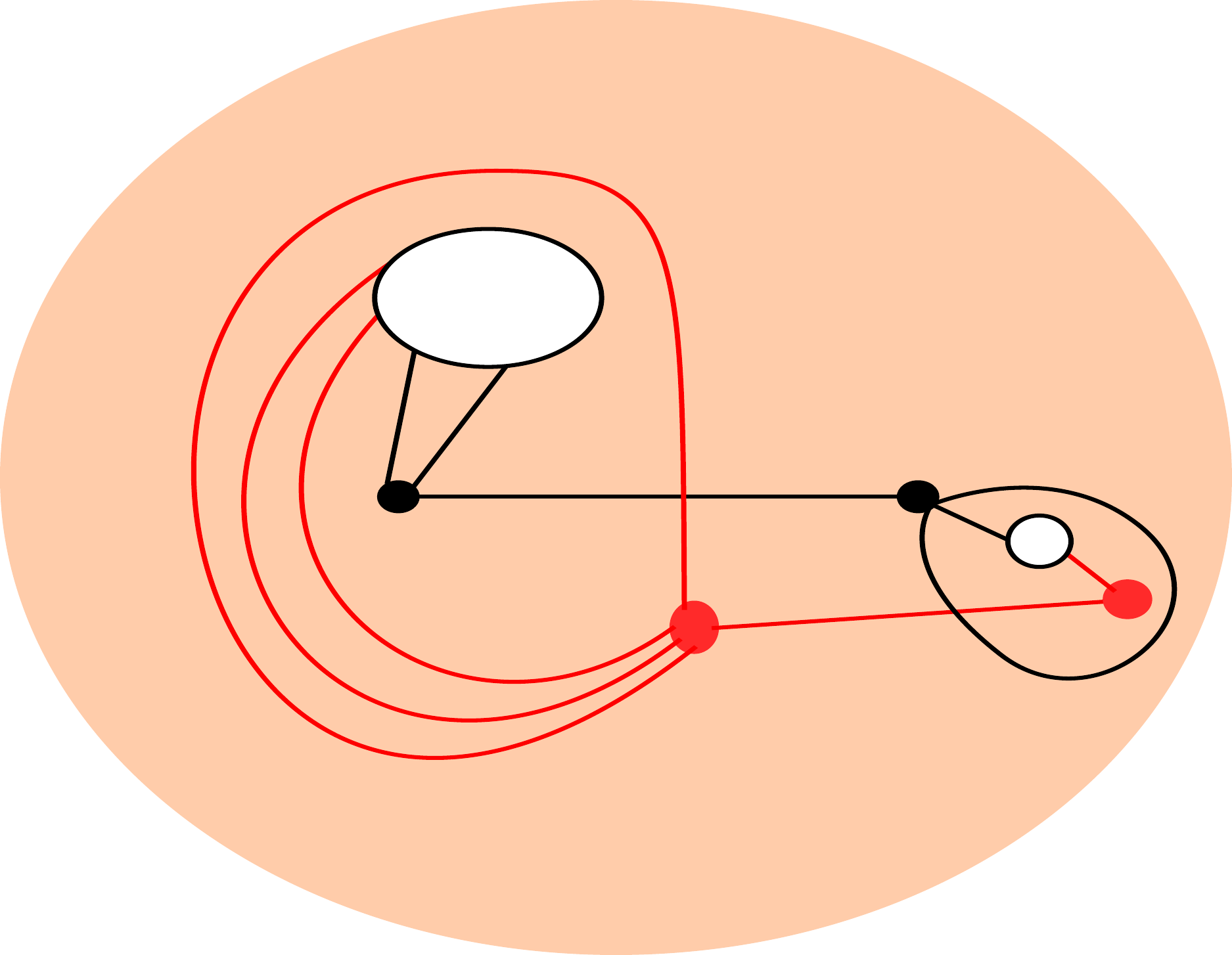}
\caption{ {\small The geometric  dual (in red) of a 
$h$-weak HEG  cellular embedding  (in black) in a punctured sphere. }}
\label{fig:dualhesdg}
\end{minipage}
\end{figure}

 Illustrations of geometric duals of HEG $h$-proper and $h$-weak cellular
embeddings are given in Figures \ref{fig:dualheg} and
\ref{fig:dualhesdg}.

The geometric duality for cellular embeddings is know to be an involution,
i.e. $(\cG^*)^*=\cG$. We investigate if this property  is 
preserved for our present duality on $h$-weak HEG cellular embeddings. One
must notice that we only need to understand how to apply twice the duality 
reflects at the level of the external cycles. Indeed, any 2-cell of $\bar \cG \subset \Sigma \cup \bsig$ 
which does not intersect any boundary circle will map to itself (up to graph isomorphism) applying
twice the duality. 

\begin{theorem}\label{theo:dual}
Let $\cG_h\subset \Sigma$ be a   $h$-weak HEG cellular embedding.
 Then, $(\cG_{h}^*)^*\subset \Sigma$ is a $h$-weak HEG cellular embedding and $(\cG_{h}^*)^*\subset \Sigma $ is equivalent to $\cG_h\subset \Sigma $.   
\end{theorem}
 \proof 
The first statement is given by Theorem \ref{Propo:dual} applied on $\cG_{h}^*\subset \Sigma$
which is a $h$-weak HEG cellular embedding. 

Call $\overline{(\cG^*)^*}$ the completed graph of $(\cG^*_h)^*$, 
and $\bar\cG$ the completed graph of $\cG_h$. 
To prove the equivalence relation between $(\cG_{h}^*)^*\subset \Sigma$ and 
$\cG_{h}\subset \Sigma$, it is sufficient to show that 
there is an isomorphism $  \bar\cG \to \overline{(\cG^*)^*}$  which will naturally extend 
to an homeomorphism. We focus on external cycles. 
Pick a 2-cell of $\Sigma - \phi(\cG)$ which contains boundary circles
where are incident some completed half-edges.  
Let us call this external cycle $c_{\ext}$, with vertex set $\{v_k\}\cup \cV_\circ$, 
with $\cV_\circ$ its set of leaves intersecting boundary circles.
By construction, we know that $c_\ext$  corresponds to a unique special vertex $v^*$ in the 
dual $\cG^*_h\subset \Sigma$ and $v^*$ is adjacent to leaves intersecting boundary circles. 
The set of those leaves which is denoted by $\cV_{\circ}^*$ is partitioned in 
the same way that $\cV_\circ$ is partitioned on boundary circles reflecting the property of $\cG_h \subset \Sigma$ to be $h$-weak. 
We need to show that applying again the duality at this vertex $v^*$ leads us
back to an external cycle which is equal to the initial one $c_{\ext}$, up to graph isomorphism. 
The edges incident to $v^*$ which are not incident to leaves 
in $\cV_{\circ}^*$ will map in $(\cG^*)^*$ to edges which will close a cycle $c'$ around
$v^*$. The edges of $c'$ are incident to vertices $\{\tilde{v}_k\}$ in a way 
that the corresponding edges of $c_{\ext}$ are incident to vertices $\{v_k\}$. 
It is clear that $\{\tilde{v}_k\}$ and $\{v_k\}$ are in one-to-one correspondence
and the edges incident to these are also one-to-one and the incidence
relation of a chain-type graph is preserved. 
The rest of the procedure becomes straightforward because
the graph at the external cycle is planar: for each $\tilde{v}_k$, we construct a cycle corresponding
to a part $s_{c_l}$ of $\cV_0^*$ associated with leaves incident to the same boundary circle $c_l$. 
Grafing this cycle $s_{c_l}$ leads to leaves incident to $\tilde{v}_k$ in an equivalent
way that a part of $V_\circ$ was incident to $v_k$. 

\qed

\section{Duality and polynomial invariants }
\label{sect:dual}

The Tutte polynomial has the fundamental property that, for the dual of a planar graph $\cG$, $T (\cG; X, Y ) = T (\cG^*; Y, X)$. Bollob\'as and Riordan derived a similar result for the BR polynomial invariant after restricting of some its variables  \cite{bollo}.  We want to investigate the analog relation for HERGs.

\begin{definition}[Internal and external half-edges and edges]
\label{def:invariant} Consider a HERG $G_H$. 

A HR $h_0$ of $G_H$ is internal if $h_0$ is the only HR in the 
external cycle of $G_H$ containing $h_0$. Otherwise, $h_0$ is called external.

A non-loop edge $e$ of $G_H$ is internal if the two HR generated by $G_H\vee e$ are both internal; $e$ is called semi-internal if one of these HRs is internal and the second is external. Otherwise, $e$ is called external. 

A vertex $v$ of $G_H$ is external if there is at least one HR incident to $v$. Otherwise it is internal. The number of internal and external vertices are denoted by $V_{\inter}$ and $V_{\ext}$, respectively.

\end{definition}
 
\begin{definition}\label{dualherg}
Let $G_H=G(V,E,H)$ be a HERG associated with the $h$-proper HEG cellular embedding $\cG_H$. 
 The dual of $G_H$, denoted $G_H^*$, is the HERG associated with the geometric dual of
 $\cG_H$.
\end{definition}

From these definitions, we can establish the following correspondences between a HERG $G_H$ and its dual $G_H^*$:
\bea
V_{\inter}(G_H)=f_{\inter}(G_H^*),\,\,\, f_{\inter}(G_H)=V_{\inter}(G_H^*),\,\,\,
V_{\ext}(G_H)=C_{\ext}(G_H^*),\,\,\, C_{\ext}(G_H)=V_{\ext}(G_H^*).
\eea
and $e(G_H)=e(G_H^*)$, and they have an equal  number of HRs.

The  dual of HERGs as stated in \cite{krf} is written with combinatorial maps with fixed points. The construction of this dual HERG coincides for several examples with
the construction of the dual HERG as stated in Definition \ref{dualherg}. We 
therefore conjecture that these definition can be shown equivalent.

As discussed in \cite{bollo}, a bridge in a ribbon graph $G$ corresponds to a trivial loop in $G^*$  and an ordinary edge in $G$ may correspond to non-trivial loop in $G^*$. 
This property remains true for HERGs. 
We note that if $e$ is a trivial twisted (respectively untwisted) loop the contraction of $e$ gives one vertex (respectively two vertices) possibly with HRs.
 As in the case of ribbon graphs, we have the following relation: 
\bea\label{eq:dualrel}
(G_H-e)^*=G_H^*/e,\,\,\, (G_H/e)^*=G_H^*-e,
\eea
where $G_H-e$ is the deletion of the edge $e$, and $G_H/e$ is its contraction.

A polynomial invariant on HERGs was introduced in \cite{avohou}:
\begin{definition}[BR polynomial for HERGs]\label{def:polycut}
Let $G_H$ be a HERG.  
We define the polynomial of $G_H$ to be
\beq
\cR_{G_H}(x,y,z,s,w,t)
=\sum_{A\sset G_H} (x-1)^{r(G_H)-r(A)}y^{n(A)}
z^{k(A)-f_{\inter}(A)+n(A)} \, s^{C_\ext(A)} w^{o(A)}\, t^{|H(A)|},
\label{brfla}
\eeq
where the sum is performed over the spanning cutting subgraphs. The quantities $r(A)$, $n(A)$,
$k(A)$, $f_{\inter}(A)$ and $C_\ext(A)$ are respectively the rank, the nullity, the number of connected 
components, the number of closed faces and external cycles of $A$. We define $o(A)=0$ is $A$ orientable, and 1 otherwise, $|H(A)|$ is the number of HRs of $A$, and where $w^2=w$ holds. 
\end{definition}

The polynomial \eqref{brfla} satisfies the contraction/cut recurrence relation where the cut operation replaces the deletion in the context of HERGs. The cut of an edge $e$ of a HERG $G_H$, denoted by 
$G_H\vee e$, is the deletion of the edge $e$ and the insertion of two HRs attached to its incident vertices or vertex in the loop situation. Cutting a subset of edges in a HERG yields a spanning cutting subgraph of $G_H$. We have for an ordinary edge $e$: 
\beq
\cR_{G_H} = \cR_{G_H\vee e }+\cR_{G_H/e}
\eeq
For special edges (loops and bridges) reduced relation exists and can be found 
in \cite{avohou}. Although this relation looks similar to the contraction/deletion of BR, 
we must emphasize that this polynomial is not an evaluation of the BR polynomial
simply because it is defined on HERGs which have more combinatorial properties
that usual ribbon graphs. In short, 
the universality theorem of the BR polynomial for  ribbon graphs does not apply
to HERGs.

To find relations between polynomial invariants evaluated at
dual HERGs, we adopt the same strategy as in \cite{bollo}. 
We need to seek pertinent restrictions or modifications of  $\cR_{G_H}$ \eqref{brfla}. 
After restrictions, it appears possible to find some relationships. 
Two interestings cases are discussed below. 

\

\noindent{\bf Polynomial of the first kind -} 
Let us denote by $R_{G_H}(.)$ the polynomial obtained from Definition \ref{def:polycut} by replacing the spanning cutting subgraphs by the spanning subgraphs. Therefore the number of HRs remains constant in all subgraphs, and if $H=\emptyset$, then $R_{G_{H}}$ reduces to the BR polynomial. Using techniques developed for HERGs in \cite{avohou}, we can
show that, for an ordinary edge $e$, 
\bea\label{eq:twovarcontrdel}
R_{G_H}=R_{G_H/e}+R_{G_H- e}\,. 
\eea

We now introduce the following two-variable polynomial
\bea\label{eq:polypdeletion}
P_{G_H}(a,b)=\sum_{A\subset G_H}a^{f_{\inter}(A)}b^{C_\ext(A)},
\eea
where the summation is over the spanning subgraphs. 

Like $R_{G_H}$, the polynomial $P_{G_H}$ obeys a contraction/deletion recursion 
on HERGs.  At this point, one may wonder if, by the universality theorem of Bollob\'as and Riordan \cite{bollo},
 $R_{G_H}$ or $P_{G_H}$ are evaluations of the BR polynomial. 
The answer of that question is no. Both polynomials are defined on HERGs,
and we can show that they fails to satisfy the vertex  union
operation on simple examples. As a consequence, all known recipe theorems worked out for ribbon graphs cannot be used here. 
In the following, this will be further explained as it 
will appear clear that the universality theorem 
for BR polynomial cannot be applied neither for $R_{G_H}$  nor for  $P_{G_H}$.

Let $E_n$ be the HERG made with $n$ isolated vertices possibly with HRs. We have $P_{E_n}(a,b)=a^{n-C_\ext(E_n)}b^{C_\ext(E_n)}=P_{E_n^*}(a,b)$ since $E_n$ is self-dual. Setting $Q_{G_H}=P_{G_H^*}$ and using \eqref{eq:dualrel}, for every edge $e$ we have:
\bea
P_{G_H^*- e}(a,b)+P_{G_H^*/e}(a,b)\,
=P_{(G_H/e)^*}(a,b)+P_{(G_H- e)^*}(a,b)\,
=Q_{G_H/e}(a,b)+Q_{G_H- e}(a,b).
\eea
Hence 
\bea\label{eq:polystabledual}
P_{G_H}(a,b)=P_{G_H^*}(a,b).
\eea

We now consider the case of one-vertex HERG. We have  $f_{\inter}(A)-1=n(A)-[k(A)-f_{\inter}(A)+n(A)]$ and 
\bea\label{eq:sumspanning}
P_{G_H}(a,b)=a\sum_{A\subset G_H}a^{f_{\inter}(A)-1}b^{C_\ext(A)}=aR_{G_H}(x,a,a^{-1},b,1,1).
\eea
From a direct calculation, one gets
\bea\label{eq:pcontrat}
P_{G_H}(a,b) = \left\{\begin{array}{ll} 
(b+1)P_{G_H/e}(a,b) &  {\text{if}}\,\, e \,\, {\text{is an external bridge}}\\\\
(a+1)P_{G_H/e}(a,b) &  {\text{otherwise}}\\
\end{array} \right. 
\eea
and 
\beq
\label{eq:rcontrat}
R_{G_H}(x,a,a^{-1},b,1,1) = \left\{\begin{array}{ll} 
\big(a^{-1}b(x-1)+1\big)R_{G_H/e}(x,a,a^{-1},b,1,1) &  {\text{if}}\,\, e \,\, {\text{is an external bridge}}\\\\
xR_{G_H/e}(x,a,a^{-1},b,1,1) &  {\text{otherwise}}\\
\end{array} \right. 
\eeq
Choosing $x-1=a$,  we observe that \eqref{eq:pcontrat} and \eqref{eq:rcontrat} coincide. 
Furthermore,  $R_{G_H}(.)$ and $P_{G_H}(.)$ satisfy  the contraction deletion recurrence relation for any ordinary edge and this leads to 
\begin{theorem}\label{Theo:duality1} Let $G_H$ be a HERG and $G_H^*$ its geometric dual. Then
\beq\label{duality1}
R_{G_H}(a+1,a,a^{-1},b,1,1)=R_{G_H^*}(a+1,a,a^{-1},b,1,1).
\eeq
\end{theorem}
Finally, we note that neither \eqref{eq:pcontrat} nor \eqref{eq:rcontrat} 
satisfy a single relation for general bridges. This implies that
 $R_{G_H}$ and $P_{G_H}$ fall out of the hypothesis of the universality theorem of BR polynomial. At $a=b$, that is in the limit where external cycles and
 closed faces are not distinguished, \eqref{eq:pcontrat} and \eqref{eq:rcontrat}
 merges into a single relation each, and then 
 \eqref{duality1} reduces to the duality 
 relation  of the BR polynomial on ribbon graphs obtained in \cite{bollo}.

\

\noindent{\bf Polynomial of the second kind -} There is another  
restriction of $\cR_{G_H}$ which could be mapped to $\cR_{G_H^*}$.

Consider the two-variable polynomial defined by  
\bea\label{eq:sumspanningcutn}
\cP_{G_H}(a,b)=\sum_{A\sset G_H}a^{f_{\inter}(A)}b^{C_\ext(A)},
\eea
considered as an element of the quotient of $\mathbb{Z}[a,b]$ by the ideal generated
by $b^2-ab$. As a consequence of this relation, if $n\geq0$ and $m>0$,  
$a^nb^m=b^{n+m}$. Still, if $m=0$, then the monomial $a^n$ might occur. 
Thus, in the case of $A=G_H$ with $H=\emptyset$, the monomial $a^{f_{\inter}(G_H)}$ 
appears in this expansion of $\cP_{G_H}$ and this term cannot be reduced. 
We write 
\bea
&&
\cP_{G_H}(a,b) = M_{G_H}(a,b)
+\sum_{A\sset \cG/ A \notin \{ G_H\vee E(G_H), G_H\}} a^{f_{\inter}(A)}b^{C_\ext(A)} \,, 
\crcr
&&
M_{G_H}(a,b) =  
\sum_{ A \in \{ G_H\vee E(G_H), G_H\}} a^{f_{\inter}(A)}b^{C_\ext(A)}\,, 
\eea
where $M_{G_H}(a,b)$ is the sum of the contributions of $G_H$
and  $G_H\vee E(G_H)$, the HERG obtained by cutting all edges in $G_H$. The monomial associated with $G_H$ is in fact the same as that of $G_H/E(G_H)$, the HERG obtained by contracting all edges in $G_H$.  

For one-vertex HERGs, we have the relation
\bea\label{eq:sumspanningcut}
\cP_{G_H}(a,b)=a\sum_{A\sset G_H}a^{f_{\inter}(A)-1}b^{C_\ext(A)}
=a\cR_{G_H}(x,a,a^{-1},b,1,1). 
\eea
 The polynomial $\cP_{G_H}$ satisfies the contraction/cut relation for any edge $e$, 
that is
\beq
\cP_{G_H} = \cP_{G_H \vee e} + \cP_{G_H /e}\,. 
\eeq
The above relation is a corollary of Theorem 3.11 of \cite{avohou}. 

\begin{lemma}\label{cpp}
Let $G_H$ a HERG with $C_\ext(G_H) >0$, then 
\beq
\cP_{{G_H}}(a,b) = \cP_{{G_H}}(b,b) = P_{{G_H}}(b, b)\,. 
\eeq
\end{lemma}
\proof Using $C_\ext(G_H) >0$, we can show that all 
monomial in the state sum of $\cP_{G_H}$ can be mapped
to $b^{f_\inter(A) + C_{\ext}(A)}$, hence the first equality. 
The second equality follows from the fact that
the spanning subgraphs of $G_H$ are in one-to-one
correspondence with spanning cutting subgraphs of $G_H$
and the fact that the quantity $ f_{\inter}(A)+ C_{\ext}(A)$ remains the same for the corresponding subgraphs. 

\qed

\begin{proposition}\label{prop:polydualp}
Let $G_H$ be a HERG.  

\begin{itemize}
\item[(1)]  If $C_{\ext}(G_H)> 0$, $M_{G_H}(a,b)=M_{G_H^*}(a,b)$ and 
$\cP_{G_H}(a,b)=\cP_{G_H^*}(a,b)$.

\item[(2)] If $C_{\ext}(G_H)= 0$, $M_{G_H}(a,b)=M_{G_H^*}(b,a)$ and $\cP_{G_H}(a,b)-M_{G_H}(a,b)=\cP_{G_H^*}(a,b)-M_{G_H^*}(a,b)$.
\end{itemize}
\end{proposition}
\proof
The resulting graph $G_H/E(G_H)$ gives some isolated vertices possibly with HRs. The vertices without HRs are in one-to-one correspondence with the internal faces of $G_H$ and the  remaining are one-to-one with the external faces of $G_H$. The contribution of $G_H/E(G_H)$ in $M_{G_H}(.)$ is 
$a^{f_{\inter}(G_H)}b^{C_\ext(G_H)}$ and the contribution of $G_H\vee E(G_H)$ is $b^{V(G_H)}$,
with $V(G_H)$ the number of vertices of $G_H$.

We start by proving (1).   Suppose that $C_{\ext}({G_H})> 0$. Using $ab=b^2$,
$M_{{G_H}}(a,b)=a^{f_{\inter}(G_H)}b^{C_\ext(G_H)}+b^{V(G_H)}=b^{C_\ext(G_H)+f_{\inter}(G_H)}+b^{V(G_H)}=b^{V(G_H^*)}+b^{C_\ext(G_H^*)+f_{\inter}(G_H^*)}=M_{{G_H^*}}(a,b)$.

Similarly, $\cP_{{G_H}}(a,b)=\sum_{A\sset {G_H}}b^{f_{\inter}(A)+C_\ext(A)}$. 
Then using Lemma \ref{cpp}, we have $\cP_{{G_H}}(a,b)$ $= P_{{G_H}}(b,b)$. The fact that $P_{{G_H}}(b,b)=P_{{G_H^*}}(b,b)$ \eqref{eq:polystabledual}
yields the result.

We now prove (2).  We now suppose that $C_{\ext}({G_H})= 0$.
Then $M_{{G_H}}(a,b)=a^{f_{\inter}({G_H})}+b^{V({G_H})}$. Knowing that the vertices of ${G_H}^*$ are in one-to-one correspondence with the faces of ${G_H}$, then 
\bea\label{eq:mm}
M_{{G_H}^*}(a,b)=a^{f_{\inter}({G_H}^*)}+b^{V({G_H}^*)}=a^{V({G_H})}+b^{f_{\inter}({G_H})}=M_{{G_H}}(b,a).
\eea
Let us prove that 
\beq\label{eq:sssetcut}
\sum_{A\sset G_H/ A \notin \{ G_H\vee E(G_H), G_H\}}b^{f_{\inter}(A)+C_\ext(A)}=
\sum_{A\sset G_H^*/ A \notin \{ G_H\vee E(G_H), G_H\}}b^{f_{\inter}(A)+C_\ext(A)}\,.
\eeq
Let us embed $P_{{G_H}}(a,b)$ in $\mathbb{Z}[a,b]/ \langle b^2-ab\rangle$. 
$C_\ext(G_H)=0$ implies that $H=\emptyset$ and $G_H=G$ becomes a ribbon graph.  Then 
\bea
P_{{G_H}}(a,b) - a^{f_{\inter}({G_H})}- a^{V({G_H})}
= \sum_{A\subset {G_H}/  A \notin \{ G_H- E(G_H), G_H\}}a^{f_{\inter}(A)+C_\ext(A)}.
\eea
where $G_H- E(G_H)$ is the ribbon graph obtained from $G_H$ after deleting all its
edges. Using \eqref{eq:polystabledual} and the equality $a^{f_{\inter}({G_H})}+a^{V({G_H})}=a^{V({G_H^*})}+a^{f_{\inter}({G_H^*})}$, we certainly have 
\bea\label{eq:subsetsset}
P_{{G_H}}(a,b) - a^{f_{\inter}({G_H})}- a^{V({G_H})}
=P_{{G_H^*}}(a,b) - a^{f_{\inter}({G_H^*})}- a^{V({G_H^*})}\,. 
\eea
We can achieve \eqref{eq:sssetcut} using \eqref{eq:subsetsset} evaluated at $a=b$ 
and Lemma \ref{cpp}. 

\qed

\

Coming back to relation \eqref{eq:sumspanningcut}, we can now explore the bridge case. We have 
\beq\label{eq:pcontratcut}
\cP_{{G_H}}(a,b) = \left\{\begin{array}{ll} 
(a^{-1}b^2+1)\cP_{{G_H}/e}(a,b) &  {\text{if}}\,\, e \,\, {\text{is an internal bridge}}\\\\
(b+1)\cP_{{G_H}/e}(a,b) &  {\text{otherwise}}\\
\end{array} \right. 
\eeq
and 
\beq\label{eq:rcontratcut}
\cR_{{G_H}}(x,a,a^{-1},b,1,1) = \left\{\begin{array}{ll} 
\big(a^{-2}b^2(x-1)+1\big)\cR_{{G_H}/e}(x,a,a^{-1},b,1,1) &  {\text{if}}\, e \,\, {\text{is an internal bridge}}\\\\
\big(a^{-1}b(x-1)+1\big)\cR_{{G_H}/e}(x,a,a^{-1},b,1,1) &  {\text{otherwise}}\\
\end{array} \right. 
\eeq
Once again choosing $x-1=a$,   \eqref{eq:pcontratcut} and \eqref{eq:rcontratcut} 
delivers the same information. The recursion relation of contraction/cut satisfied by  $\cR_{{G_H}}(.)$ and $\cP_{{G_H}}(.)$ for any ordinary edge 
and Proposition \ref{prop:polydualp}  lead us to the following statement. 

\begin{theorem}\label{theo:dualerg}
Let ${G_H}$ be a HERG. 
\begin{itemize}
\item[(1)] If $C_{\ext}({G_H})> 0$, then
$$\cR_{{G_H}}(a+1,a,a^{-1},b,1,1)=\cR_{{G_H}^\star}(a+1,a,a^{-1},b,1,1).$$

\item[(2)] If $C_{\ext}({G_H})= 0$, then
$$a\cR_{{G_H}}(a+1,a,a^{-1},b,1,1)-M_{{G_H}}(a,b)=a\cR_{{G_H}^\star}(a+1,a,a^{-1},b,1,1)-M_{{G_H}^*}(a,b).$$
\end{itemize}
\end{theorem}

\section*{Acknowledgments}
R.C.A acknowledges the support of the Trimestre Combinatoire ``Combinatorics and Interactions",
Institut Henri Poincar\'e, Paris, France, at initial stage of this work. 
R.C.A is partially supported by the Third World Academy of Sciences (TWAS)
and the Deutsche Forschungsgemeinschaft (DFG) through a TWAS-DFG Cooperation Visit
grant. Max-Planck Institute for Gravitational Physics, Potsdam, Germany,
is thankfully acknowledged for its hospitality. 
The ICMPA is in partnership  with  the  Daniel  Iagolnitzer  Foundation  (DIF), France.

\vspace{0.5cm}

%\begin{center}
%\rule{3cm}{0.01cm}
%\end{center}

\end{document}